\documentclass{IEEEtran}  % FINAL VERSION (12pp)

\IEEEoverridecommandlockouts                              % This command is only
\def\TITLE{Distributed reactive power feedback control for voltage regulation and loss minimization}
\def\TITLEnl{Distributed reactive power feedback control\\for voltage regulation and loss minimization}
\def\AUTHORS{Saverio Bolognani, Ruggero Carli, Guido Cavraro, and Sandro Zampieri}

\usepackage{amsmath}
\usepackage{amsfonts}
\usepackage{amssymb}
\usepackage{amsthm}

\newtheorem{theorem}{Theorem}
\newtheorem{corollary}{Corollary}
\newtheorem{lemma}{Lemma}
\newtheorem{proposition}{Proposition}
\newtheorem{remark}{Remark}
\newtheorem{definition}{Definition}
\newtheorem{assumption}{Assumption}

\def \graph	{\mathcal{G}}		% grafo elettrico
\def \nodes	{\mathcal{V}}		% insieme dei nodi
\def \compensators {{\mathcal{C}}}		% insieme dei nodi generatori
\def \edges	{\mathcal{E}}		% insieme delle coppie di nodi comunicanti
\newcommand{\neighbors}[1]{\mathcal{N}(#1)}	% insieme dei vicini (si possono ottenere info)

\def \nonodes {n}
\def \nocompensators {m}
\def \noedges {{|\edges|}}

\def \complexnumbers {\mathbb{C}}
\def \realnumbers {\mathbb{R}}

\DeclareMathOperator{\diag}{diag}

\def\1{{\mathbf{1}}}

\usepackage[english]{babel}

\usepackage[final]{graphicx} 
\usepackage[export]{adjustbox}

\usepackage[hyperindex=true, %
  pdftitle={\TITLE},%
  pdfauthor={\AUTHORS},%
  colorlinks=true,%
  pagebackref=false,%
  plainpages=false,%
  pdfpagelabels,%
  linkcolor=black,%
  citecolor=black,%
  filecolor=black,%
  urlcolor=black%
]{hyperref}

\title{\LARGE \bf \TITLEnl}

\author{\AUTHORS%
\thanks{S. Bolognani is with the Laboratory for Information and Decision Systems, Massachusetts Institute of Technology, Cambridge, MA 02139 USA. Email: {\tt\small saverio@mit.edu}}% <-this % stops a space
\thanks{G. Cavraro, R. Carli, and S. Zampieri are with the Department of Information Engineering, University of Padova, Italy. 
Email: {\tt\small \{carlirug, guido.cavraro, zampi\}@dei.unipd.it}.} %
}

\begin{document}

\maketitle

%%%%%%%%%%%%%%%%%%%%%%%%%%%%%%%%%%%%%%%%%%%%%%%%%%%%%%%%%%%%%%%%%%%%%%%%%%%%%%%%

\begin{abstract}
We consider the problem of exploiting the microgenerators dispersed in the power distribution network
in order to provide distributed reactive power compensation for power losses minimization and voltage regulation.
In the proposed strategy, microgenerators are smart agents that can measure their phasorial voltage, share these data with the other agents on a cyber layer, and adjust the amount of reactive power injected into the grid,
according to a feedback control law that descends from duality-based methods 
applied to the optimal reactive power flow problem.
Convergence to the configuration of minimum losses and feasible voltages is proved analytically 
for both a synchronous and an asynchronous version of the algorithm, where agents update their state independently one from the other.
Simulations are provided in order to illustrate the performance and the robustness of the algorithm, and the innovative feedback nature of such strategy is discussed.
\end{abstract}

%%%%%%%%%%%%%%%%%%%%%%%%%%%%%%%%%%%%%%%%%%%%%%%%%%%%%%%%%%%%%%%%%%%%%%%%%%%%%%%%

\section{Introduction}

Recent technological advances, together with environmental and economic challenges, have been motivating the deployment of small power generators in the low voltage and medium voltage power distribution grid.
The availability of a large number of these generators in the distribution grid can yield relevant benefits to the network operation, which go beyond the availability of clean, inexpensive electrical power.
They can be used to provide a number of ancillary services that are of great interest for the management of the grid \cite{Katiraei_2006_Powermanagementstrategies,Prodanovic_2007_HarmonicandReactive}.

We focus in particular on the problem of optimal reactive power compensation for power losses minimization and voltage regulation.
In order to properly command the operation of these devices, the distribution network operator is required to solve an \emph{optimal reactive power flow} (ORPF) problem.
Powerful solvers have been designed for the ORPF problem, 
and advanced optimization techniques have been recently specialized for this task \cite{Zhao_2005_multiagent-basedparticleswarm,Villacci2006,Lavaei_2011_Powerflowoptimization}.
However, this approach assumes that an accurate model of the grid is available, that all the grid buses are monitored, that loads announce their demand profiles in advance, and that
generators and actuators can be dispatched on a day-ahead, hour-ahead, and real-time basis.
For this reason, these solvers are in general offline and centralized, and they collect all the necessary field data, compute the optimal configuration, and dispatch the reactive power production at the generators.

These tools cannot be applied directly to the ORPF problem faced in low/medium voltage power distribution networks.
The main reasons are that not all the buses of the grid are monitored, individual loads are unlikely to announce they demand profile in advance,
the availability of small size generators is hard to predict (being often correlated with the availability of renewable energy sources).
Moreover, the grid parameters, and sometimes even the topology of the grid, are only partially known,
and generators are expected to connect and disconnect, requiring an automatic reconfiguration of the grid control infrastructure (the so called \emph{plug and play} approach).

Different strategies have been recently proposed in order to address these issues. Purely local algorithms have been proposed, in which each generator
is operated according to its own measurements \cite{Prodanovic_2007_HarmonicandReactive}, in order to compensate for the voltage rise caused by its own active power injection \cite{Carvalho2008,Keane2011}. Because of the absence of coordination between microgenerators, the full potential of the microgenerators for voltage regulation is not exploited in these strategies \cite{Vovos2007}.

Different coordination strategies have been then proposed, for example by casting the problem into the framework of resource allocation \cite{Baran2007} and by using hierarchical dispatch schemes \cite{Rogers2010}. A two-stage approach has been proposed in \cite{Robbins2013}, where microgenerators first attempt to regulate their voltage autonomously, and they involve their neighbors in this task if their regulation capability is not sufficient.

Finally, some distributed approaches that do not require any central controller, but still require measurements at all the buses of the distribution grid, have been proposed. In order to derive a distributed algorithm for this problem, different convex relaxation methods \cite{Jabr2008,Bai2008,Zhang2013a,Lavaei2014} have been applied, and various distributed optimization algorithms have been specialized for the resulting convex ORPF problem \cite{Lam2012,Farivar2012,DallAnese2013,Sulc2013}.

Only recently, algorithms that are truly scalable in the number of generators and do not require the monitoring of all the buses of the grid, have been proposed for the problem of power loss minimization (with no voltage constraints) \cite{Tenti2012,Bolognani2013,Bolognani2013w}.
While these algorithms have been designed by specializing classical nonlinear optimization algorithms to the ORPF problem, they can also be considered as \emph{feedback} control strategies.
Indeed, the key feature of these algorithms is that they require the alternation of measurement and actuation based on the measured data, and therefore they are inherently online algorithms.
In particular, the reactive power injection of the generators is adjusted by these algorithms based on the phasorial voltage measurements that are performed at the buses where the generators are connected.
The resulting closed loop system features a tight dynamic interconnection of the \emph{physical layer} (the grid, the generators, the loads) with the \emph{cyber layer} (where communication, computation, and decision happen).
In this paper, we design a distributed feedback algorithm for the ORPF problem with voltage constraints, in which the goal is to minimize reactive power flows while ensuring that the voltage magnitude across the network lies inside a given interval, and that the microgenerators' reactive power limits are not violated. The analysis of the convergence of this strategy is based on an assumption of homogeneity of the X/R ratio of the power lines across the network. The robustness of the proposed solution with respect to possible variability in these parameters has been investigated via simulations.

In Section \ref{sec:problem_formulation}, a cyber-physical model for a smart power distribution grid is provided.
In Section \ref{sec:ORPF}, the ORPF problem with voltage constraints is formulated.
A feedback control strategy for its solution is derived in Section~\ref{sec:dual}, by using the tools of dual decomposition.
A synchronous and an asynchronous version of the algorithm are presented in Section~\ref{sec:dual}  and Section~\ref{sec:asynchronous}, respectively.
The convergence of both the proposed algorithms is studied in Section~\ref{sec:convergence}.
Some simulations are provided in Section~\ref{sec:simulations}, while Section~\ref{sec:conclusion} concludes the paper discussing some relevant features of the feedback nature of the proposed strategy.

%%%%%%%%%%%%%%%%%%%%%%%%%%%%%%%%%%%%%%%%%%%%%%%%%%%%%%%%%%%%%%%%%%%%%%%%%%%%%%%%
\section{Mathematical preliminaries and notation}
\label{sec:notation}

Let $\graph= (\nodes, \edges, \sigma,\tau)$ be a directed graph, 
where $\nodes$ is the set of nodes,
$\edges$ is the set of edges, 
and $\sigma, \tau: \edges \rightarrow \nodes$ are two functions such that
edge $e \in \edges$ goes from the source node $\sigma(e)$ to the terminal node $\tau(e)$. 

Given two nodes $h,k \in \nodes$, we define the path $\mathcal P_{hk}$ as the sequence of adjacent nodes, without repetitions, that connect node $h$ to node $k$.

Let $A \in \{0, \pm 1\}^{\noedges \times \nonodes}$ 
be the incidence matrix of the graph $\graph$, defined via its elements
$$
[A]_{ev} = \left\{\begin{array}{cl}
-1 & \text{if }v=\sigma(e) \\
1 & \text{if }v=\tau(e) \\
0 & \text{otherwise.}
\end{array}\right.
$$
If the graph $\graph$ is connected 
(i.e. for every pair of nodes there is a path connecting them),
then $\1$ is the only vector in the null space $\ker A$,
$\1$ being the column vector of all ones.
We define by $\1_v$ the vector whose value is $1$ in position $v$, and $0$ everywhere else.

In the rest of the paper we will often introduce complex-valued functions defined on nodes and on edges.
These functions will also be intended as vectors in $\complexnumbers^\nonodes$ (where $\nonodes = |\nodes|$) and $\complexnumbers^\noedges$.
Given a vector $u$, we denote by $\bar{u}$ its (element-wise) complex conjugate, and by $u^T$ its transpose. 
We denote by $\Re(u)$ and $\Im(u)$ the real and imaginary part of $u$, respectively.

%%%%%%%%%%%%%%%%%%%%%%%%%%%%%%%%%%%%%%%%%%%%%%%%%%%%%%%%%%%%%%%%%%%%%%%%%%%%%%%%

\section{Cyber-physical model of a smart power distribution grid}
\label{sec:problem_formulation}

In this work, we envision a \emph{smart} power distribution network as a cyber-physical system, in which
\begin{itemize}
\item the \textbf{physical layer} consists of the power distribution infrastructure, including power lines, loads, microgenerators, and the point of connection to the transmission grid;
\item the \textbf{cyber layer} consists of intelligent agents, dispersed in the grid, and provided with actuation, sensing, communication, and computational capabilities.
\end{itemize}

\begin{figure}[bt]
\centering
\resizebox{74mm}{!}{
\includegraphics[width=85mm]{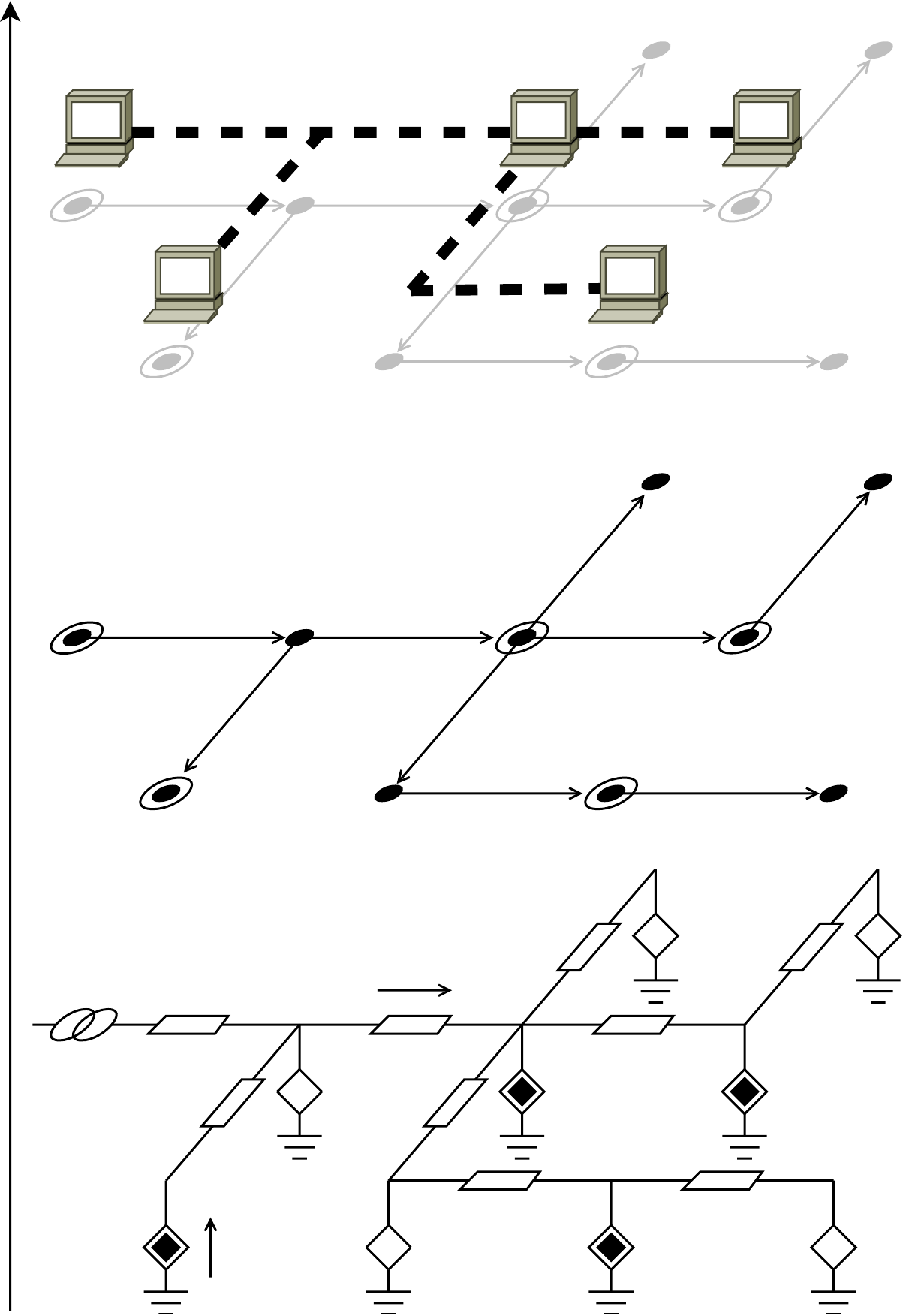}
\put(-250,25){\rotatebox{90}{\small Physical layer}}
\put(-250,150){\rotatebox{90}{\small Graph model}}
\put(-250,265){\rotatebox{90}{\small Cyber layer}}
\put(-213,35){$u_v$}
\put(-181,14){$i_v$}
%\put(-136,94){$\xi_e$}
\put(-138,65){$z_e$}
\put(-214,138){$v$}
\put(-170,190){$\sigma(e)$}
\put(-118,190){$\tau(e)$}
\put(-140,170){$e$}
\put(-224,190){$0$}
\put(-229,333){agent}
\put(-189,323){communication}
}
\caption{Schematic representation of the power distribution grid model.
In the lower panel the \emph{physical layer} is represented via a circuit representation, where black diamonds are microgenerators, white diamonds are loads, and the left-most element of the circuit represents the PCC.
The middle panel illustrates the adopted \emph{graph representation} for the same grid. Circled nodes represent both microgenerators and the PCC.
The upper panel represents the \emph{cyber layer}, where agents (i.e. microgenerator nodes and the PCC)  are also connected via some communication infrastructure.
}
\label{fig:microgrid_model}
\end{figure}
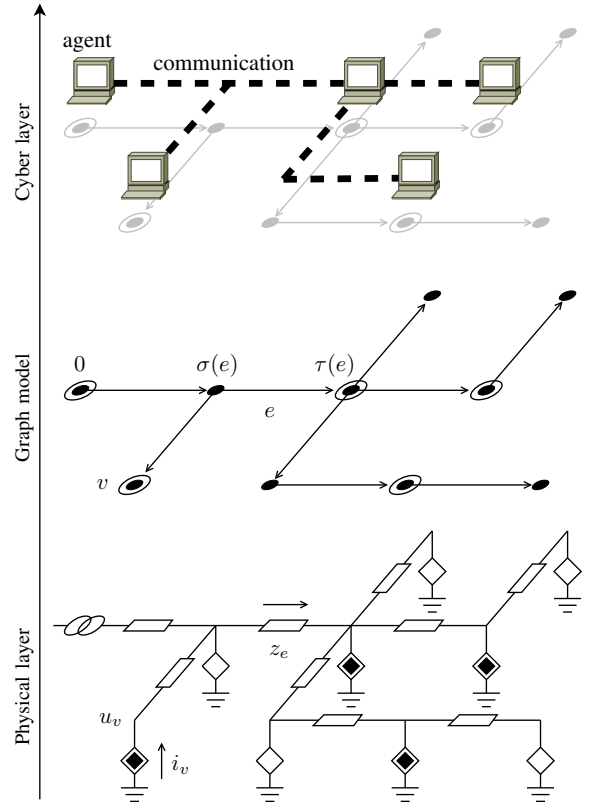

\subsection{Physical layer}

For the purpose of this paper, we model the physical layer of a smart power distribution network as a directed graph $\graph$, 
in which edges represent the power lines, and nodes represent buses (see Figure~\ref{fig:microgrid_model}, middle panel).
Buses correspond to loads, microgenerators, and also the point of connection to the transmission grid (called point of common coupling, or PCC, and indexed as node $0$).

We limit our study to the steady state behavior of the system, when all voltages and currents are sinusoidal signals at the same frequency.
Each signal can therefore be represented via a complex number $y = |y|e^{j\angle y}$ whose absolute value $|y|$ corresponds to the signal root-mean-square value, 
and whose phase $\angle y$ corresponds to the phase of the signal with respect to an arbitrary global reference.
In this notation, the steady state of the grid is described by the following system variables (see Figure~\ref{fig:microgrid_model}, lower panel):
\begin{itemize}
\item $u \in \complexnumbers^\nonodes$, where $u_v$ is the grid voltage at node $v$;
\item $i \in \complexnumbers^\nonodes$, where $i_v$ is the current injected at node $v$.
\end{itemize}

We model the grid power lines as series impedances, neglecting their shunt admittance. For every edge $e$ of the graph, we define by $z_e$ the impedance of the corresponding power line.
We assume the following.
\begin{assumption}
\label{ass:theta}
All the power lines in the grid have the same inductance/resistance (X/R) ratio, but possibly different impedance magnitude, i.e.
$$
z_e = e^{j\theta} |z_e|
$$
for any $e$ in $\edges$ and for a fixed $\theta$.
\end{assumption}

This assumption is satisfied when the X/R ratio of the power lines of the grid is relatively homogeneous, which is reasonable in many practical cases (see for example the IEEE standard testbeds \cite{Kersting2001}).  In Section~\ref{sec:simulations} we investigate what is the effect of a possible variability of the X/R ratio, and we show how the proposed strategy is very robust against this possible uncertainty.

Under this assumption, we have a linear relation between
bus voltages and currents in the form
\begin{equation}
i = e^{-j\theta} L u
\label{eq:iLu}
\end{equation}
where $L:= A^T Z^{-1} A$ is the weighted Laplacian of the graph, in which $A$ is the incidence matrix of  $\graph$, and $Z = \diag(|z_e|, e\in \edges)$ is the diagonal matrix of the magnitudes of line impedances.

Each node $v$ of the grid is then characterized by a law 
relating its injected current $i_v$ with its voltage $u_v$.
We model the PCC as an ideal sinusoidal voltage generator at the microgrid nominal voltage $U_N$ 
with arbitrary fixed angle $\psi$
\begin{equation}
u_0 = U_N e^{j \psi}.
\label{eq:PCCmodel}
\end{equation}
In the power system analysis terminology, node $0$ is then a \emph{slack bus} with fixed voltage magnitude and angle.

We model loads and microgenerators (that is, every node $v$ of the microgrid except the PCC) via the following law relating the voltage $u_v$ and the current $i_v$
\begin{equation}
u_v\bar i_v = s_v, %
%\left|\frac{u_v}{U_N}\right|^{\eta_v},
\quad \forall v \in \nodes \backslash \{ 0 \},
\label{eq:ZIPModel}
\end{equation}
where $s_v$ is the injected \emph{complex power}. 
The quantities
$$
p_v := \Re(s_v) \quad \text{and} \quad q_v := \Im(s_v) 
$$
are denoted as \emph{active} and \emph{reactive} power, respectively.
The complex powers $s_v$ corresponding to grid loads are such that $\{p_v <0\}$, meaning that positive active power is \emph{supplied} to the devices.
The complex powers corresponding to microgenerators, on the other hand, are such that $\{p_v \ge 0\}$, as positive active power is \emph{injected} into the grid.
In the power system analysis terminology, all nodes but the PCC are being modeled as \emph{constant power} or \emph{P-Q buses}.
Microgenerators fit in this model, as they generally are commanded via a complex power reference and they can inject it independently from the voltage at their point of connection
\cite{Lopes_2006_Definingcontrolstrategies,Green_2007_Controlofinverterbased}.

\subsection{Cyber layer}

\begin{figure}[bt]
\centering
\includegraphics[trim=22mm 0mm 3mm 0mm,clip,scale=0.85,valign=t]{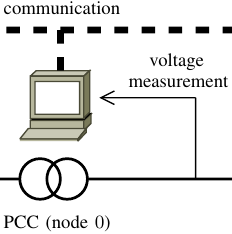}\hspace*{\stretch{1}} 
\includegraphics[trim=3mm 0mm 3mm 0mm,clip,scale=0.85,valign=t]{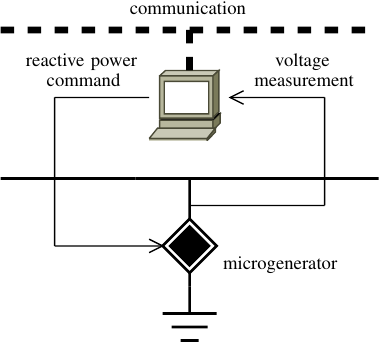}
\caption{A schematic representation of the agents' capabilities and of the way in which agents of the cyber layer interface with the physical layers. The first panel represent the agent at the PCC (node $0$), which is provided with voltage measurement capabilities. The second panel represents all the other agents, at the microgenerators, which are provided with both measurement and actuation capabilities. All the agent can communicate via some communication channel.}
\label{fig:agents}
\end{figure}

We assume that every microgenerator, and also the PCC, correspond to an \emph{agent} in the cyber layer (upper panel of Figure~\ref{fig:microgrid_model}).
We denote by $\compensators$ (with $|\compensators|=m$) this subset of the nodes of $\graph$.
Each agent is provided with some computational capability, and with some sensing capability, in the form of 
a phasor measurement unit (i.e. a sensor that can measure voltage amplitude and angle \cite{Phadke_1993_Synchronizedphasormeasurements}).
Agents that corresponds to a microgenerator can also actuate the system, by commanding a set point for the amount of reactive power injected by that microgenerator
(see Figure~\ref{fig:agents}).

Finally, agents can communicate, via some communication channel that could possibly be the same power lines
(via power line communication).
Motivated by this possibility, we define the neighbors in the cyber layer in the following way.

\begin{figure}[tb]
\centering
\footnotesize
\includegraphics[width=6.1cm]{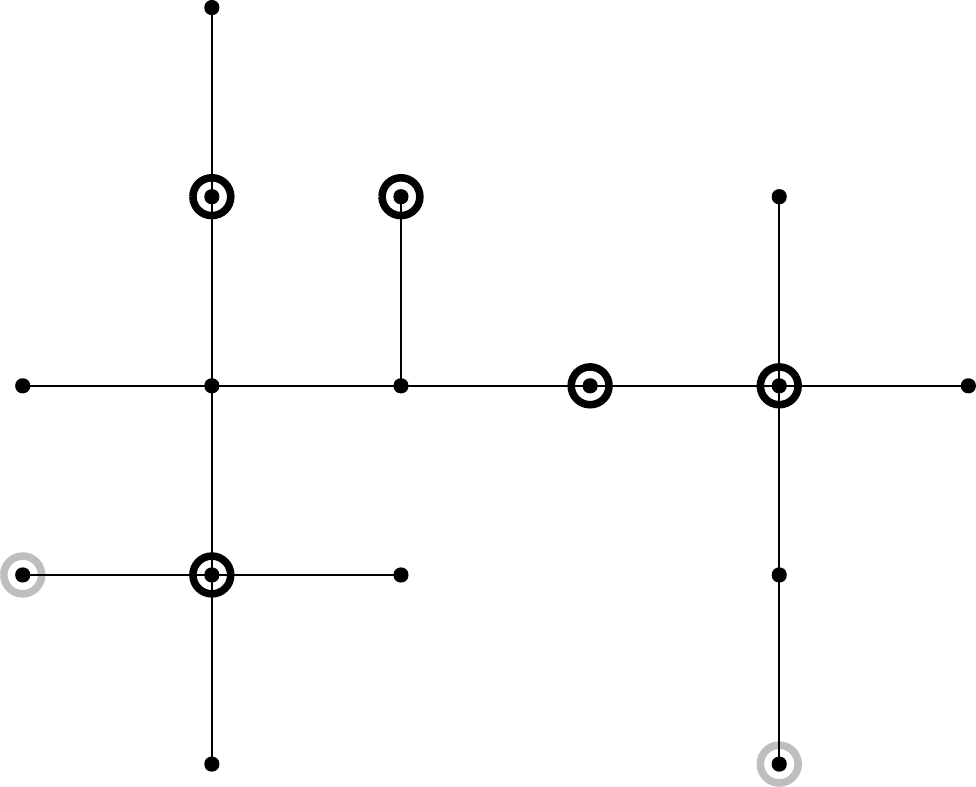}
\put(-71,81){$h$}
\put(-120,116){$k \in \neighbors{h}$}
\put(-195,47){$k' \notin \neighbors{h}$}
\caption{An example of neighbor agents in the cyber layer. Circled nodes (both gray and black) are agents (nodes in $\compensators$). Nodes circled in black belong to the set $\neighbors{h} \subset \compensators$.
Node circled in gray are agents which do not belong to the set of neighbors of $h$. 
For each agent $k \in \neighbors{h}$, the path that connects $h$ to $k$ does not include any other agent besides $h$ and $k$ themselves.}
\label{fig:neighbors}
\end{figure}

\begin{definition}[Neighbors in the cyber layer]
Let $h \in \compensators$ be an agent of the cyber layer. The set of agents that are neighbors of $h$, denoted as $\neighbors{h}$, is the subset of $\compensators$ defined as
$$
\neighbors{h} = \left\{ k \in \compensators \;|\; \forall \mathcal P_{hk}, \mathcal P_{hk} \cap \compensators = \{h,k\}\right\}.
$$
\label{neighbors}
\end{definition}
Figure~\ref{fig:neighbors} gives an example of such set.
We assume that every agent $h\in \compensators$ knows its set of neighbors $\neighbors{h}$, and can communicate with them.
Notice that this architecture can be constructed by each agent in a distributed way, 
for example by exploiting the PLC channel (as suggested for example in \cite{Costabeber_2011_Ranging}).
This allows also a plug-and-play reconfiguration of such architecture when new agents are connected to the grid.

\section{Approximated model and its properties}

In this section we review an approximate explicit solution of the nonlinear equations \eqref{eq:iLu}, \eqref{eq:PCCmodel}, and \eqref{eq:ZIPModel} which has been proposed in \cite{Bolognani2013}. This approximation will play a crucial role in deriving our distributed control strategy to solve the optimal reactive power flow problem that we will introduce in next section.
In order to present the approximated solution, we need the following technical lemma.

\begin{lemma}[Lemma 1 in \cite{Bolognani2013}]
Let $L$ be the weighted Laplacian of $\graph$.
There exists a unique symmetric, positive semidefinite matrix $X \in \realnumbers^{n\times n}$ such that
\begin{equation}
\begin{cases}
XL = I - \1 \1_0^T \\
X \1_0 = 0.
\end{cases}
\label{eq:Xproperties}
\end{equation}
\label{lemma:X}
\end{lemma}

The matrix $X$ depends only on the topology of the grid power lines and on their impedance.
The matrix $X$ has some notable properties, including the fact that 
\begin{equation}
X_{hh} \ge X_{hk} \ge 0 \quad h,k \in \nodes,
\label{eq:signX}
\end{equation}
and the fact that 
\begin{equation}
(\1_h - \1_k)^T X (\1_h - \1_k) = \left|Z^\text{eff}_{hk} \right|, \quad h,k \in \nodes,
\label{eq:zeff}
\end{equation}
where $Z^\text{eff}_{hk}$ represents the \emph{effective impedance} of the power lines between node $h$ and $k$.
Notice that, if the grid is radial (i.e. $\graph$ is a tree) then $Z^\text{eff}_{hk}$ is simply the impedance of the only path from node $h$ to node $k$.
Now let us introduce the following block decomposition of the vector of voltages $u$
$$
u = \begin{bmatrix}
u_0 \\ u_G \\ u_L
\end{bmatrix},
$$
where $u_0$ is the voltage at the PCC,
$u_G \in \complexnumbers^{m-1}$ are the voltages at the microgenerators, 
and $u_L \in \complexnumbers^{n-m}$ are the voltages at the loads.
Similarly, we also define $s_G = p_G + j q_G$ and $s_L = p_L + j q_L$.

By adopting this block decomposition as before, we have
\begin{equation}
X = \begin{bmatrix}
0 & 0 & 0\\
0 & M & N \\
0 & N^T & Q 
\end{bmatrix},
\label{eq:Xblocks}
\end{equation}
with $M \in \realnumbers^{(m-1)\times(m-1)}$, $N \in \realnumbers^{(m-1)\times(m-n)}$, and $Q \in \realnumbers^{(n-m)\times(n-m)}$.
The following proposition provides the approximate relation between the grid voltages and the power injections at the nodes.

\begin{proposition}
Consider the physical model described by the set of nonlinear equations %\eqref{eq:KCL}, \eqref{eq:KVL}, 
\eqref{eq:iLu}, \eqref{eq:PCCmodel}, and \eqref{eq:ZIPModel}.
Node voltages then satisfy
\begin{equation}
\begin{bmatrix}
u_0 \\ u_G \\ u_L
\end{bmatrix}
\!\!=
e^{j\psi}
\!\!
\left( 
U_N \1 + 
\frac{e^{j\theta}}{U_N}
\!\!
\begin{bmatrix}
0 & 0 & 0 \\ 
0 & M & N \\
0 & N^T & Q
\end{bmatrix}
\!\!
\begin{bmatrix}
0 \\ 
\bar s_G \\
\bar s_L
\end{bmatrix}
\right)
\!
+
o\left(\frac{1}{U_N}\right)
\label{eq:approximate_solution}
\end{equation}
where the little-o notation means that $\lim_{U_N\rightarrow \infty} \frac{o(f(U_N))}{f(U_N)} = 0$.
\label{pro:approximation}
\end{proposition}
\begin{IEEEproof}
It descends directly from Proposition 1 in \cite{Bolognani2013}.
\end{IEEEproof}

The quality of this approximation relies on having large nominal voltage $U_N$ and relatively small currents injected by the inverters (or supplied to the loads).
This assumption is verified in practice, and corresponds to correct design and operation of power distribution networks,
where indeed the nominal voltage is chosen sufficiently large 
(subject to other functional constraints) 
in order to deliver electric power to the loads with relatively small power losses on the power lines.
The model proposed in Proposition~\ref{pro:approximation} extends the DC power flow model \cite[Chapter 3]{Gomez_2009_Electricenergysystems} to the case in which lines are not purely inductive. This way, the model is able to describe the voltage drop on the lines, and therefore also the corresponding power losses, in a form that is conveniently linear in the complex power injections and demands.
We conclude this section by introducing the following matrix $G$.

\begin{lemma}
There exists a unique symmetric matrix $G \in \realnumbers^{m\times m}$ such that
$$
\begin{cases}
\begin{bmatrix}
0 & 0 \\ 0 & M
\end{bmatrix} G = I - \1 \1_0^T \\
G \1 = 0.
\end{cases}
$$
\label{lem:G}
\end{lemma}
\begin{IEEEproof}
The following matrix $G$ satisfies the conditions.
\begin{equation}
G = \begin{bmatrix}
 \1^T M^{-1}\1 & -\1^T M^{-1} \\ -M^{-1}\1 & M^{-1}
\end{bmatrix}.
\label{eq:G}
\end{equation}
The proof of uniqueness, that we omit here, follows exactly the same steps as in the proof of Lemma \ref{lemma:X}.
\end{IEEEproof}
The matrix $G$ has also a remarkable sparsity pattern, as the following lemma states.
\begin{lemma}\label{lem:sparsity}
The matrix $G$ has the sparsity pattern induced by the Definition \ref{neighbors} of neighbor agents in the cyber layer, i.e.
$$
G_{hk} \ne 0 \quad \Leftrightarrow \quad k \in \neighbors{h}.
$$
\end{lemma}
The proof is provided in the Appendix \ref{app:G}, where we also discuss how the elements of $G$ can be estimated by the agents, given a local knowledge of the power grid parameters.

\section{Optimal reactive power flow problem}
\label{sec:ORPF}

We consider the problem of commanding the reactive power injection of the microgenerators in order to minimize power distribution losses on the power lines
and to guarantee that the voltage magnitude and the reactive power injection stay within pre-assigned intervals.
The decision variables (or, equivalently, the inputs of the system) are therefore the reactive power setpoints $q_h, h \in \compensators\backslash\{0\}$, compactly written as $q_G$.

Power distribution losses can be expressed as a function of the voltage drop on the lines and therefore, in a matricial quadratic form, as
%, by using \eqref{eq:KVL}, as
\begin{equation}
\quad J_\text{losses} := 
%\sum_{e\in \edges} |\xi_e|^2 \Re(z_e) = 
\bar u^T L u .
\label{eq:exactlosses}
\end{equation}
Given a lower bound $U_\text{min}$ and an upper bound $U_\text{max}$ for the voltage magnitudes, and a lower bound $q_\text{min}$ and an upper bound $q_\text{max}$ for the reactive power injected by each microgenerator, we can therefore formulate the following optimization problem,
\begin{subequations}
\begin{align}
%\min_{q_h, h \in \compensators\backslash\{0\}} & \quad \bar u^T L u
\min_{q_G} & \quad \bar u^T L u
\label{eq:ORPF1} \\
\text{subject to} & \quad 
\begin{array}{l}
|u_h| \ge U_\text{min}\\
|u_h| \le U_\text{max}\\
q_h \ge q_\text{min}\\
q_h \le q_\text{max}
\end{array}
\quad \forall h \in \compensators\backslash\{0\}\label{eq:ORPF2}%
\end{align}
\label{eq:ORPF}%
\end{subequations}%
where voltages $u$ are a function of the decision variables $q_G,$ 
via the implicit relation defined by the system of nonlinear equations \eqref{eq:iLu}, \eqref{eq:PCCmodel}, and \eqref{eq:ZIPModel}.
From a control design prospective, the system-wide problem that we are considering is therefore characterized by 
\begin{itemize}
\item the input variables $q_G$,
\item the measured output variables $\begin{bmatrix} u_0 \\ u_G \end{bmatrix}$,
\item the unmeasured disturbances $p_L$, $q_L$, $p_G$.
\end{itemize}

The goal of this paper is to design a control algorithm to tackle the ORPF problem in a distributed fashion, where each microgenerator $h$ is allowed to communicate only with its neighbors in the cyber layer, i.e., the agents in $\mathcal{N}(h)$.

\begin{remark}
While the decision variables of the ORPF problem (i.e. the input variables $q_G$) do not include the reactive power supplied by the PCC (i.e. $q_0 = u_0 \bar i_0$), this quantity will also change every time the reactive power setpoints of the generators are updated by the algorithm,
because the inherent physical behavior of the slack bus (the PCC) ensures that equations \eqref{eq:iLu}, \eqref{eq:PCCmodel} and \eqref{eq:ZIPModel} are satisfied at every time.
\end{remark}

\begin{remark}
In the above formulation we have assumed that all the microgenerators have to satisfy the same reactive power injection constraint. This scenario can be seamlessly extended to the case of heterogeneous microgenerators, where $q_\text{h,min}\le q_h \le q_\text{h,max}$ being the values $q_\text{h,min}$, $q_\text{h,max}$, $h \in \compensators$, in general different for different microgenerators. 
\end{remark}

\section{A synchronous algorithm based on dual decomposition}
\label{sec:dual}

In this section, in order to design a distributed feedback control strategy to solve the ORPF problem, we apply the tool of dual decomposition to \eqref{eq:ORPF}.
Specifically, we use the approximate explicit solution of the nonlinear equations \eqref{eq:iLu}, \eqref{eq:PCCmodel}, and \eqref{eq:ZIPModel} introduced in Proposition~\ref{pro:approximation}, to derive update steps for a \emph{dual ascent algorithm} \cite{Bertsekas1999} that can be implemented distributively by the agents and that can be used as a feedback control update law.

It is convenient to reformulate the problem \eqref{eq:ORPF} via a change of coordinates, obtaining
\begin{subequations}
\begin{align}
%\min_{q_h, h \in \compensators\backslash\{0\}} & \quad \bar u^T L u
\min_{q_G} & \quad \bar u^T L u
 \label{eq:ORPFF1} \\
\text{subject to} & \quad 
\begin{array}{l}
v_h \ge v_\text{min}\\
v_h \le v_\text{max}\\
w_h \ge w_\text{min}\\
w_h \le w_\text{max}
\end{array}
\quad \forall h \in \compensators\backslash\{0\}\label{eq:ORPFF2}
\end{align}
\label{eq:ORPFF}
\end{subequations}
where 
\begin{align*}
v_h &= {|u_h|^2}/{U_N^2} & w_h &= {2 q_h}/{U_N^2} \\
v_\text{min} &= {U_\text{min}^2}/{U_N^2} &
	w_\text{min} &= {2 q_\text{min}}/{U_N^2} \\
v_\text{max} &= {U_\text{max}^2}/{U_N^2} & 
	w_\text{max} &= {2 q_\text{max}}/{U_N^2}.
\end{align*}
Basically, in \eqref{eq:ORPFF2} with the respect to \eqref{eq:ORPF2},
we have squared and normalized the constraints on the voltage magnitude and we have normalized the constraints on the power injection. While these modifications does not have any effect on the optimization problem, they will allow us to simplify the derivation of the algorithm we are going to present. 
The Lagrangian of the problem \eqref{eq:ORPFF} is
\begin{align}
&\mathcal L(q_G,\lambda_ \text{min}, \lambda_ \text{max}, \mu_ \text{min}, \mu_ \text{max})  \nonumber\\
& =\bar u^T L u + \lambda_\text{min}^T \left(v_\text{min} \1 - v_G \right)+ \lambda_{\text{max}}^T \left(v_G - v_\text{max} \1 \right) \nonumber \\
&\qquad \,\,\,\,\,\,\,+\mu_\text{min}^T \left(w_\text{min} \1 - w_G \right)+ \mu_{\text{max}}^T \left(w_G - w_\text{max} \1 \right)
\label{equ:lagrangian}
\end{align}
where $\lambda_\text{min}$, $\lambda_\text{max}$, $\mu_\text{min}$, $\mu_\text{max}$ are the Lagrangian multipliers (i.e. the dual variables of the problem) and $u$, $v_G$, $w_G$ are functions of the decision variables $q_G$, even if the dependence has been omitted. To have a more compact notation let 
$$
\nu=\left[\lambda^T_\text{min}\,\,\,\lambda^T_\text{max} \,\,\,\mu^T_\text{min}\,\,\, \mu^T_\text{max} \right]^T.
$$

A dual ascent algorithm consists in the iterative execution of the following alternated steps
\begin{enumerate}
\item dual gradient ascent step on the dual variables
\begin{align*}
&\lambda_\text{min}(t+1) = \left[\lambda_\text{min}(t) + \gamma \frac{\partial \mathcal L(q_G(t),\nu(t))}{\partial \lambda_\text{min}} \right]_+\\
&\lambda_\text{max}(t+1) = \left[\lambda_\text{max}(t) + \gamma \frac{\partial \mathcal L(q_G(t),\nu(t))}{\partial \lambda_\text{max}} \right]_+\\
&\mu_\text{min}(t+1) = \left[\mu_\text{min}(t)  + \gamma \frac{\partial \mathcal L(q_G(t),\nu(t))}{\partial \mu_\text{min}} \right]_+\\
&\mu_\text{max}(t+1) = \left[\mu_\text{max}(t)+ \gamma \frac{\partial \mathcal L(q_G(t),\nu(t))}{\partial \mu_\text{max}} \right]_+,
\end{align*}
where the $[\cdot]_+$ operator corresponds to the projection on the positive orthant,
and where $\gamma$ is a suitable positive constant;
\item minimization of the Lagrangian with respect to the primal variables $q_G$
\begin{align}\label{eq:stepPrimal}
& q_G(t+1) = \arg\min_{q_G} \mathcal L(q_G,\nu(t+1)).
\end{align}
\end{enumerate}

Observe that the updates of the Lagrange multipliers can be performed naturally in a distributed way by the agents, based on their local measurement of the violation of the voltage and power constraints. Indeed let $\lambda_{\text{min},h}$, $\lambda_{\text{max},h}$, $\mu_{\text{min},h}$ and $\mu_{\text{max},h}$ be the components of the the Lagrange multipliers $\lambda_\text{min}$, $\lambda_\text{max}$, $\mu_\text{min}$ and $\mu_\text{max}$, respectively, related to the compensator $h$. Then it easily follows that the dual step can be implemented as
\begin{align*}
\lambda_{\text{min},h}(t+1)&=\left[\lambda_{\text{min},h}(t)+\gamma(v_\text{min}-v_h)  \right]_+\\
\lambda_{\text{max},h}(t+1)&=\left[\lambda_{\text{max},h}(t)+\gamma(v_h-v_\text{max})  \right]_+\\
\mu_{\text{min},h}(t+1)&=\left[\mu_{\text{min},h}(t)+\gamma(w_\text{min}-w_h)  \right]_+\\
\mu_{\text{max},h}(t+1)&=\left[\mu_{\text{max},h}(t)+\gamma(w_h-w_\text{max})  \right]_+.
\end{align*}
The crucial point is to derive an expression for the minimizer $q_G(t+1)$ in \eqref{eq:stepPrimal} that can be computed distributively by the compensators. To do so, we exploit the approximation introduced in Proposition~\ref{pro:approximation} obtaining a value of $q_G(t+1)$ which  is equivalent to the one in \eqref{eq:stepPrimal} up to a term which vanishes to zero for large nominal voltage $U_N$.  

The resulting update corresponds to the algorithm that we now present.
We assume here that the agents are coordinated, i.e., they can update their state variables $q_h$, $\lambda_{\text{min},h}$, $\lambda_{\text{max},h}$, $\mu_{\text{min},h}$ and $\mu_{\text{max},h}$, synchronously.

\noindent\hrulefill\\
\noindent\textbf{Synchronous algorithm}

Let all agents (except the PCC) store the auxiliary scalar variables $\lambda_{\text{min},h}$, $\lambda_{\text{max},h}$, $\mu_{\text{min},h}$ and $\mu_{\text{max},h}$.
Let $\gamma$ be a positive scalar parameter,
and let $\theta$ be the impedance angle defined in Assumption~\ref{ass:theta}.
Let $G_{hk}$ be the elements of the sparse matrix $G$ defined in Lemma~\ref{lem:G}.
At every synchronous iteration of the algorithm, each agent $h \in \compensators\backslash\{0\}$ executes the following operations in order:
\begin{itemize}
\item	it measures its voltage $u_h$ and it gathers the voltage measurements $$\{u_k = |u_k| \exp(j\angle u_k), k \in \neighbors{h}\}$$ from its neighbors;
\item 	it updates the auxiliary variables $\lambda_{\text{min},h}$, $\lambda_{\text{max},h}$, $\mu_{\text{min},h}$, $\mu_{\text{max},h}$, as
\begin{align*}
& \lambda_{\text{min},h} \leftarrow \left[\lambda_{\text{min},h} + \gamma \left(\frac{U_\text{min}^2}{U_N^2} - \frac{|u_h|^2}{U_N^2}\right)\right]_+\\
& \lambda_{\text{max},h} \leftarrow \left[\lambda_{\text{max},h} + \gamma \left(\frac{|u_h|^2}{U_N^2}- \frac{U_\text{max}^2}{U_N^2}\right)\right]_+\\
& \mu_{\text{min},h} \leftarrow \left[\mu_{\text{min},h} + \gamma \left(\frac{q_\text{min}}{U_N^2} - \frac{q_h}{U_N^2}\right)\right]_+\\
& \mu_{\text{max},h} \leftarrow \left[\mu_{\text{max},h} + \gamma \left(\frac{q_h}{U_N^2}- \frac{q_\text{max}}{U_N^2}\right)\right]_+;
%\label{eq:lambda-sync}
\end{align*}
\item it gathers from its neighbors the updated values of the Lagrange multipliers $\mu_{\text{min},k}$, $\mu_{\text{max},k}$, $k \in \neighbors{h}$;
\item based on the new values of $\lambda_{\text{max},h}$, $\lambda_{\text{min},h}$ and of $\mu_{\text{min},k}$, $\mu_{\text{max},k}$, $k \in \neighbors{h}$, it updates the injected reactive power $q_h$ as
\begin{equation}
\begin{split}
q_h \leftarrow \ &q_h - \sin \theta (\lambda_{\text{max},h}  - \lambda_{\text{min},h}) \\
				&\,\,\,+ \sum_{k \in \neighbors{h}} G_{hk} |u_h| |u_k| \sin(\angle u_k - \angle u_h - \theta)\\
				&\,\,\,- \sum_{k \in \neighbors{h} \backslash \{0\}} G_{hk} (\mu_{\text{max},k}  - \mu_{\text{min},k}) .
\end{split}
\label{eq:q-sync}
\end{equation}
\end{itemize}
\hrulefill

Observe that the above algorithm can be implemented in a completely distributed fashion. Indeed each agent is required to exchange information only with its neighbors in the cyber layer. 

The following Proposition shows how the update \eqref{eq:q-sync} approximates the primal step \eqref{eq:stepPrimal}.

\begin{proposition}
Consider the synchronous algorithm above described. 
Then
\begin{align*}
\frac{\partial \mathcal L(q_G(t+1), \nu(t+1))}{\partial q_G} = o\left(\frac{1}{U_N^2}\right),
\end{align*}
namely, the update \eqref{eq:q-sync}, minimizes the Lagrangian with respect to the primal variables, up to a term that vanishes for large $U_N$.
\label{prop:lagrangian_gradient}
\end{proposition}

The details of the proof of Proposition~\ref{prop:lagrangian_gradient} are postponed to Appendix~\ref{app:primal}. It is based on the following technical lemma, which will be useful again later.

\begin{lemma}
Consider the Lagrangian  $\mathcal L(q_G,\lambda)$ defined in \eqref{equ:lagrangian}.
The partial derivative with respect to the primal variables $q_G$ is
\begin{equation*}
\begin{split}
\frac{\partial \mathcal L(q_G,\nu)}{\partial q_G} &=  \frac{2}{U_N^2} \big(M q_G + N q_L + \sin \theta M (\lambda_\text{max}  - \lambda_\text{min}) \\
				& \phantom{=} + \mu_\text{max}  - \mu_\text{min} \big) + o\left(\frac{1}{U_N^2}\right).
\end{split}
\end{equation*}
\label{lem:primalderivative}
\end{lemma}

Also the proof of Lemma~\ref{lem:primalderivative} is in Appendix~\ref{app:primal}.

Via these steps, we therefore specialized the dual ascent steps to the ORPF problem that we are considering, and we obtained a distributed feedback control law for the system. We study the convergence of the closed loop system in Section~\ref{sec:convergence}.

\begin{remark}
\label{rem:feedback}
It is important to notice that the proposed synchronous algorithm requires that the agents actuate the system at every iteration, by updating the set point for the amount of reactive power injected by the microgenerators.
Only by doing so, the subsequent measurement of the voltages will be informative of the new state of the system.
The resulting control strategy is thus a \emph{feedback strategy} that necessarily requires the real-time interaction of the controller (the cyber layer) with the plant (the physical layer), as depicted in Figure~\ref{fig:feedback}.
This tight interaction between the cyber layer and the physical layer is the fundamental feature of the proposed approach, and allows to
drive the system towards the optimal configuration, which in principle depends on the reactive power demands of the loads, without collecting this information from them.
In a sense, the algorithm is inferring this hidden information from the measurement performed on the system during its execution.

\begin{figure}[!t]
\centering
\includegraphics[width=85mm]{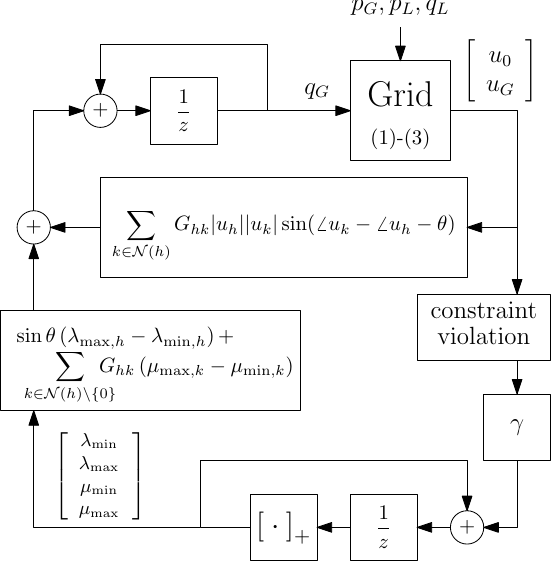}
\caption{A block diagram representation of the synchronous control algorithm proposed in Section~\ref{sec:dual}, where the tight interconnection of the cyber and the physical layer (i.e. the feedback strategy) is evident.}
\label{fig:feedback}
\end{figure}
\end{remark}

%%%%%%%%%%%%%%%%%%%%%%%%%%%%%%%%%%%%%%%%%%%%%%%%%%%%%%%%%%%%%%%%%%%%%%%%%%%%%%%%%

\section{Asynchronous algorithm}
\label{sec:asynchronous}
In order to avoid the burden of system-wide coordination among the agents, we also propose an asynchronous version of the algorithm, 
in which the agents corresponding to the microgenerators update their state $(q_h, \lambda_{\text{max},h}, \lambda_{\text{min},h}, \mu_{\text{max},h}, \mu_{\text{min},h})$
independently one from the other, based on the information that they can gather from their neighbors.

We assume that each agent (except for the agent located at the PCC) is provided with an individual timer, by which it is triggered, and that no coordination is present between these timers: they tick randomly, with exponentially, identically distributed waiting times.

\noindent\hrulefill\\
\noindent\textbf{Asynchronous algorithm}
\nopagebreak

Let all agents (except the PCC) store four auxiliary scalar variables $\lambda_{\text{max},h}, \lambda_{\text{min},h}, \mu_{\text{max},h}, \mu_{\text{min},h}$.
Let $\gamma$ be a positive scalar parameter,
and let $\theta$ be the impedance angle defined in Assumption~\ref{ass:theta}.
Let $G_{hk}$ be the elements of the matrix $G$ defined in Lemma~\ref{lem:G}.

When agent $h \in \compensators\backslash\{0\}$ is triggered by its own timer, it performs the following actions in order:
\begin{itemize}
\item	it measures its voltage $u_h$ and it gathers from its neighbors
the voltage measurements 
$$\{u_k = |u_k| \exp(j\angle u_k), k \in \neighbors{h}\}$$ 
and the values of the Lagrange multipliers
$$\{\mu_{\text{min},k}, \mu_{\text{max},k},\, k \in \neighbors{h}\};$$ 
\item 	it updates the auxiliary variables $\lambda_{\text{min},h}$, $\lambda_{\text{max},h}$, $\mu_{\text{min},h}$, $\mu_{\text{max},h}$, as
\begin{align*}
& \lambda_{\text{min},h} \leftarrow \left[\lambda_{\text{min},h} + \gamma \left(\frac{U_\text{min}^2}{U_N^2} - \frac{|u_h|^2}{U_N^2}\right)\right]_+\\
& \lambda_{\text{max},h} \leftarrow \left[\lambda_{\text{max},h} + \gamma \left(\frac{|u_h|^2}{U_N^2}- \frac{U_\text{max}^2}{U_N^2}\right)\right]_+\\
& \mu_{\text{min},h} \leftarrow \left[\mu_{\text{min},h} + \gamma \left(\frac{q_\text{min}}{U_N^2} - \frac{q_h}{U_N^2}\right)\right]_+\\
& \mu_{\text{max},h} \leftarrow \left[\lambda_{\text{max},h} + \gamma \left(\frac{q_h}{U_N^2}- \frac{q_\text{max}}{U_N^2}\right)\right]_+;
%\label{eq:lambda-sync}
\end{align*}
\item	based on the new value of $\lambda_{\text{min},h}$, $\lambda_{\text{max},h}$, it updates the injected reactive power $q_h$ as
\begin{equation}
\begin{split}
q_h \leftarrow \ &q_h - \sin \theta (\lambda_{\text{max},h}  - \lambda_{\text{min},h}) +\\
				&+ \sum_{k \in \neighbors{h}} G_{hk} |u_h| |u_k| \sin(\angle u_k - \angle u_h - \theta) +\\
				&- \sum_{k \in \neighbors{h} \backslash \{0\}} G_{hk} (\mu_{\text{max},k}  - \mu_{\text{min},k}) .
\end{split}
\end{equation}
\end{itemize}
\hrulefill

The update equations for the asynchronous algorithm are exactly the same of the synchronous case.
Here, however, we update both the primal and the dual variable of the agents independently and asynchronously.
Also the analysis of the convergence of this algorithm is postponed to the next section.

%%%%%%%%%%%%%%%%%%%%%%%%%%%%%%%%%%%%%%%%%%%%%%%%%%%%%%%%%%%%%%%%%%%%%%%%%%%%%%%%%

\section{Convergence analysis}
\label{sec:convergence}

In this section, we investigate the convergence of both the synchronous algorithm proposed in Section \ref{sec:dual} and of the asynchronous algorithm proposed in Section \ref{sec:asynchronous}. 

In order to do so, we rewrite the terms that appeared in the dual ascent update step, namely
$$
\frac{\partial \mathcal L(q_G,\nu)}{\partial \lambda_\text{min}}, \quad 
\frac{\partial \mathcal L(q_G,\nu)}{\partial \lambda_\text{max}}, \quad
\frac{\partial \mathcal L(q_G,\nu)}{\partial \mu_\text{min}}, \quad
\frac{\partial \mathcal L(q_G,\nu)}{\partial \mu_\text{max}},
$$
using the expression introduced in Proposition~\ref{pro:approximation} for the voltages. We start from 
\begin{equation*}
\frac{\partial \mathcal L(q_G,\nu)}{\partial \lambda_\text{min}}
=
v_\text{min} \1 -  v_G. 
\end{equation*}
By plugging in the approximate solution \eqref{eq:approximate_solution}, via some algebraic manipulations, we can  express $v_G$ as
\begin{align}
v_G&=\1 + \frac{2}{U_N^2} \Re\left(
e^{j\theta}
M \bar s_G 
+
e^{j\theta}
N \bar s_L 
\right)
+ o\left(\frac{1}{U_N^2}\right) \nonumber \\
&=\1 +\frac{2}{U_N^2}(\sin \theta M q_G  + \cos \theta M p_G ) \nonumber \\
&\quad +\frac{2}{U_N^2} \left(  \cos \theta N p_L  + \sin \theta N q_L \right)+ o\left(\frac{1}{U_N^2}\right) \label{eq:vg}
\end{align}
and, in turn, we have that
\begin{equation}
\frac{\partial \mathcal L(q_G,\nu)}{\partial \lambda_\text{min}}
= \frac{2}{U_N^2} \left( b_\text{min} - \sin \theta M q_G \right)+ o\left(\frac{1}{U_N^2}\right)
\end{equation}
where 
$$
b_\text{min} = \frac{U_N^2}{2} (v_\text{min} - 1)\1  - (M p_G \cos \theta  + N (p_L \cos \theta + q_L \sin \theta)).
$$
Similar calculations lead to

\begin{equation*}
\frac{\partial \mathcal L(q_G,\nu)}{\partial \lambda_\text{max}}
= \frac{2}{U_N^2} \left(\sin \theta M q_G  - b_\text{max}  \right)+ o\left(\frac{1}{U_N^2}\right),
\end{equation*}
where
$$
b_\text{max} = \frac{U_N^2}{2} (v_\text{max} - 1)\1  - (M p_G \cos \theta  + N (p_L \cos \theta + q_L \sin \theta)).
$$
Additionally observe that
\begin{align*}
\frac{\partial \mathcal L(q_G,\nu)}{\partial \mu_\text{min}}
&=
w_\text{min} \1 - w_G  = \frac{2}{U_N^2} \left(\1 q_\text{min}  - q_G  \right),\\
\frac{\partial \mathcal L(q_G,\nu)}{\partial \mu_\text{max}}
&=
w_G - w_\text{max} \1  = \frac{2}{U_N^2} \left(q_G - \1 q_\text{max} \right).
\end{align*} 
The proposed dual ascent step can therefore be rewritten in compact form as
\begin{equation}
\nu(t+1) = \left[\nu(t) + \gamma \frac{2}{U_N^2} \left(
\Phi q_G(t) + b
\right) + o\left(\frac{1}{U_N^2}\right) \right]_+,
\label{eq:updatelambda_compact}
\end{equation}

where 
\begin{equation}\label{eq:Phi}
\Phi = \begin{bmatrix}
-\sin\theta M\\\sin\theta M\\-I\\I
\end{bmatrix}, \quad b =  \begin{bmatrix}
b_\text{min}\\-b_\text{max}\\ \1q_\text{min}\\ -\1q_\text{max}
\end{bmatrix}.
\end{equation}

The update step for the primal variables can be rewritten based on Proposition~\ref{prop:lagrangian_gradient} and Lemma~\ref{lem:primalderivative}, obtaining
\begin{align}\label{eq:primal_update_compact}
q_G(t+1)& =  - M^{-1}Nq_L - M^{-1} \Phi \nu(t+1)  + o\left(\frac{1}{U_N}\right).
\end{align}

In the analysis that follows, we study the approximated description of the closed loop system in which we neglect the infinitesimal terms.
Notice that, by doing so, both the voltages $u$ and the squared voltage magnitudes $v_G$ become affine functions of the decision variables $q_G$.
By plugging those expressions in the formulation of the ORPF problem \eqref{eq:ORPFF}, one obtains the following strictly convex quadratic problem with linear inequality constraints
\begin{subequations}
\begin{align}
\min_{q_G} & \quad q_G^T \frac M 2 q_G + q_G^T N q_L \\
\text{subject to} & \quad \Phi q_G  + b \le 0,
\end{align}
\label{eq:ORPF_appx}
\end{subequations}
for which strong duality holds. The rest of the section is split into two subsections: in the first one we consider the synchronous version of the algorithm, in the second one the asynchronous version. 

\subsection{Syncronous case}\label{subsec:syn}

For the synchronous version of the algorithm, we consider the update equations
\begin{equation}
\nu(t+1) = \left[\nu(t) + \gamma \frac{2}{U_N^2} \left( \Phi q_G(t)  + b \right)\right]_+,
\label{eq:appl}
\end{equation}
for the dual variables, and
\begin{equation}
q_G(t+1) = - M^{-1}Nq_L -M^{-1} \Phi^T \nu(t+1).
\label{eq:appq}
\end{equation}
for the primal variables. Observe that \eqref{eq:appl} and \eqref{eq:appq} 
differ from \eqref{eq:updatelambda_compact} and from \eqref{eq:primal_update_compact} only by infinitesimal terms, and they correspond to the 
standard equation for the dual ascent steps for \eqref{eq:ORPF_appx}. 
Indeed, the equilibrium $(q_G^*, \nu^*)$ of \eqref{eq:appl}-\eqref{eq:appq} is characterized by
$$
\Phi q_G^*  + b  \le 0
\quad \text{and} \quad
q_G^* + M^{-1}Nq_L + M^{-1} \Phi^T \nu^*= 0,
$$
which correspond to the  necessary conditions for the optimality according to Uzawa's saddle point theorem \cite{Uzawa1958}. 

It will be useful in the following to define $\sigma_\text{min}$ and $\sigma_\text{max}$ as the minimum and the maximum eigenvalue of $M$, respectively.
The following result characterizes the convergence of the algorithm described by \eqref{eq:appl} and \eqref{eq:appq}.

\begin{theorem}
\label{th:s}
Consider the optimization problem \eqref{eq:ORPF_appx} and the dynamic system described by the update equations \eqref{eq:appl} and \eqref{eq:appq}.
Then the trajectory $t \to q(t)$ converges to the optimal primal solution $q_G^*$ if 
$$
\gamma \le \frac{U_N^2}{\rho (\Phi M^{-1} \Phi^T)},
$$
where 
\begin{align*}
&\rho(\Phi M^{-1}\Phi^T)=
 2\max\{\sigma_\text{min}^{-1} +  \sin^2 \theta \sigma_\text{min},
 \sigma_\text{max}^{-1} +  \sin^2 \theta \sigma_\text{max} \}.
\end{align*}
\end{theorem}
The proof is presented in Appendix \ref{app:convergence}.

We conclude this subsection by specializing the above result to the case where either only voltage constraints or only power constraints are considered. Observe that if we take into account only voltage constraints then the matrix $\Phi$ and the vector $b$ become
\begin{equation}\label{eq:PhiVolt}
\Phi = \begin{bmatrix}
-\sin\theta M\\\sin\theta M
\end{bmatrix}, \quad b =  \begin{bmatrix}
b_\text{min}\\-b_\text{max}
\end{bmatrix}
\end{equation}
and only the multipliers $\lambda_\text{min}$ and $\lambda_\text{max}$ are employed in the algorithm,
while if we consider only power constraints then
\begin{equation}\label{eq:PhiPower}
\Phi = \begin{bmatrix}
-I\\I
\end{bmatrix}, \quad b =  \begin{bmatrix}
 \1q_\text{min}\\ -\1q_\text{max}
\end{bmatrix}
\end{equation}
and only the multipliers $\mu_\text{min}$ and $\mu_\text{max}$ are needed. The following results follow from Theorem \ref{th:s}.
\begin{corollary}\label{cor:SyncVolt}
Consider the optimization problem \eqref{eq:ORPF_appx}, where $\Phi$ and $b$ are given as in \eqref{eq:PhiVolt}, and 
the dynamic system described by the update equations \eqref{eq:appl} and \eqref{eq:appq}.
Then the trajectory $t \to q(t)$ converges to the optimal primal solution $q_G^*$ if 
$$
\gamma \le \frac{U_N^2}{2\, \sin^2\theta\,\, \sigma_\text{max}},
$$
\end{corollary}
\begin{corollary}\label{cor:SyncPower}
Consider the optimization problem \eqref{eq:ORPF_appx}, where $\Phi$ and $b$ are given as in \eqref{eq:PhiPower}, and 
the dynamic system described by the update equations \eqref{eq:appl} and \eqref{eq:appq}.
Then the trajectory $t \to q(t)$ converges to the optimal primal solution $q_G^*$ if 
$$
\gamma \le \frac{\sigma_\text{min}\, U_N^2}{2}.
$$
\end{corollary}

\subsection{Asynchronous case}\label{subsec:asyn}

We introduce the following assumption.
\begin{assumption}
Let $\{T^{(h)}_i\}$, $i \in \mathbb N$, be the time instants in which the agent $h$ is triggered by its own timer. We assume that the timer ticks with exponentially distributed waiting times, identically distributed for all the agents in $\compensators\backslash\{0\}$.
\label{ass:timers}
\end{assumption}
Let us define the random sequence $h(t) \in \compensators\backslash\{0\}$ which tells which agent has been triggered at iteration $t$ of the algorithm.
Because of Assumption \ref{ass:timers}, the random process $h(t)$ is an i.i.d. uniform process on the alphabet $\compensators\backslash\{0\}$.
If we repeat the same analysis, neglecting the infinitesimal terms, we obtain the following update equations for the primal and dual variables, instead of 
\eqref{eq:appl} and \eqref{eq:appq}. In these equations, only the component $h(t)$ of the vectors $\lambda_\text{min}$, $\lambda_\text{max}$, $\mu_\text{min}$, $\mu_\text{max}$ and $q_G$ is updated at time $t$, namely,
\begin{equation}
\begin{split}
\lambda_{\text{min},h(t)}(t+1) &= \big[\lambda_{\text{min},h(t)}(t) + \\
&\quad +\gamma \frac{2}{U_N^2} \1_{h(t)}^T (b_\text{min} - \sin \theta M q_G(t)) \big]_+ \\
\lambda_{\text{max},h(t)}(t+1) &= \big[\lambda_{\text{max},h(t)}(t) + \\
&\quad  +\gamma \frac{2}{U_N^2} \1_{h(t)}^T (\sin \theta M q_G(t) - b_\text{max}) \big]_+ \\
\mu_{\text{min},h(t)}(t+1) &= \big[\mu_{\text{min},h(t)}(t) + \gamma \frac{2}{U_N^2}  ( q_\text{min}  - q_h(t) ) \big]_+ \\
\mu_{\text{max},h(t)}(t+1) &= \big[\mu_{\text{max},h(t)}(t) + \gamma \frac{2}{U_N^2}  ( q_h(t)-q_\text{max} ) \big]_+ 
\end{split}
\label{eq:applAsyn_h}
\end{equation}
while
\begin{equation}
\begin{split}
\lambda_{\text{min},k}(t+1) &= \lambda_{\text{min},k}(t) \\
\lambda_{\text{max},k}(t+1) &= \lambda_{\text{max},k}(t) \\
\mu_{\text{min},k}(t+1) &= \mu_{\text{min},k}(t)\\
\mu_{\text{max},k}(t+1) &= \mu_{\text{max},k}(t) 
\end{split} \qquad
\forall k\ne h(t),
\label{eq:applAsyn_noth}
\end{equation}
and 
\begin{equation}
\begin{split}
q_{h(t)}(t+1) &= - \1_{h(t)}^T ( M^{-1}Nq_L + M^{-1} \Phi^T \nu(t+1)) \\
q_k(t+1) &= q_k(t) \qquad  \forall k\ne h(t).
\end{split}
\label{eq:appqAsyn}
\end{equation}
Notice that, also in the asynchronous case, Uzawa's necessary conditions for optimality are satisfied at the equilibrium of \eqref{eq:applAsyn_h}, \eqref{eq:applAsyn_noth}  and \eqref{eq:appqAsyn}. For the asynchronous version of the algorithm we can provide theoretical results only when $\Phi$ and $b$ assume the form in \eqref{eq:PhiVolt} or in \eqref{eq:PhiPower}, namely, when we consider either only voltage constraints or only power constraints. However in the numerical section we show the effectiveness of the asynchronous algorithm when $\Phi$ and $b$ assume the general form in \eqref{eq:Phi}.

The following convergence results hold.

\begin{proposition}
\label{prop:aVolt}
Consider the optimization problem \eqref{eq:ORPF_appx}, where $\Phi$ and $b$ are given as in \eqref{eq:PhiVolt}, and the dynamic system described by the update equations \eqref{eq:applAsyn_h} and \eqref{eq:applAsyn_noth} (for the multipliers $\lambda_\text{min}$ and $\lambda_\text{max}$) and \eqref{eq:appqAsyn}. 
Let Assumption \ref{ass:timers} hold.
Then the evolution $t\to q(t)$ converges almost surely to the optimal primal solution $q_G^*$ if 
$$
\gamma \le \frac{U_N^2}{2\, \sin^2\theta\,\, \sigma_\text{max}}.
$$
\end{proposition}
\begin{proposition}
\label{prop:aPower}
Consider the optimization problem \eqref{eq:ORPF_appx}, where $\Phi$ and $b$ are given as in \eqref{eq:PhiPower}, and the dynamic system described by the update equations \eqref{eq:applAsyn_h} and \eqref{eq:applAsyn_noth} (for the multipliers $\mu_\text{min}$ and $\mu_\text{max}$) and \eqref{eq:appqAsyn}. 
Let Assumption \ref{ass:timers} hold.
Then the evolution $t\to q(t)$ converges almost surely to the optimal primal solution $q_G^*$ if 
$$
\gamma \le \frac{\sigma_\text{min}\, U_N^2}{2}.
$$
\end{proposition}

The proof is presented in Appendix \ref{app:convergence}.

%%%%%%%%%%%%%%%%%%%%%%%%%%%%%%%%%%%%%%%%%%%%%%%%%%%%%%%%%%%%%%%%%%%%%%%%%%%%%%%%%
\section{Simulations}
\label{sec:simulations}

\begin{figure}[tb]
\centering	
\includegraphics[width=0.4\textwidth]{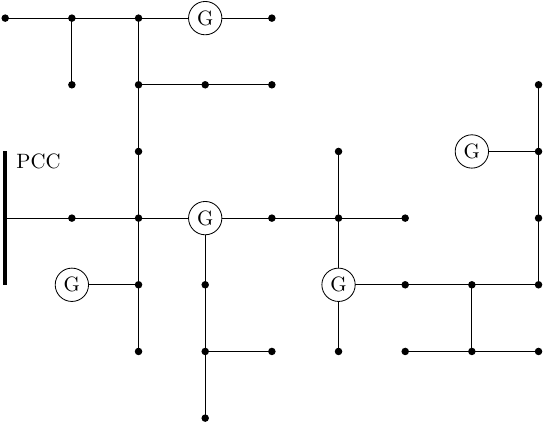}
\caption{Schematic representation of the IEEE 37 test feeder \cite{Kersting2001}, where $5$ microgenerators have been deployed.}\label{fig:ieee37}
\end{figure}

The algorithm has been tested on the testbed IEEE 37 \cite{Kersting2001},
which is an actual portion of 4.8kV power distribution network located in California.
The load buses are a blend of constant-power, constant-current, and constant-impedance loads,
with a total power demand of almost 2 MW of active power and 1 MVAR of reactive power
(see \cite{Kersting2001} for the testbed data).
The length of the power lines range from a minimum of 25 meters to a maximum of almost 600 meters.
The impedance of the power lines differs from edge to edge 
(for example, resistance ranges from 0.182 $\Omega$/km to 1.305 $\Omega$/km).
However, the inductance/resistance ratio exhibits a smaller variation,
ranging from $X/R = 0.5$ to $0.67$.
This justifies Assumption \ref{ass:theta}, in which we claimed that $\angle z_e$ can be considered constant across the network.
We considered the scenario in which $5$ microgenerators have been deployed in this portion of the power distribution grid (see Figure~\ref{fig:ieee37}).

The lower bound for voltage magnitudes has been set to 4700 V.
Both the synchronous and the asynchronous algorithm presented in Section~\ref{sec:dual} and \ref{sec:asynchronous} have been simulated on a nonlinear exact solver of the grid \cite{Zimmerman2011}.
The approximate model presented in Proposition~\ref{pro:approximation} has not been used in these simulations,
being only a tool for the design of the algorithm and for the study of the algorithm's convergence.

A time-varying  profile for the loads has been generated, in order to simulate the effect of slowly varying loads (e.g. the aggregate demand of a residential neighborhood),
fast changing / intermittent demands (e.g. some industrial loads).

\begin{figure*}[!p]
\centering	
\includegraphics[width=0.8\textwidth]{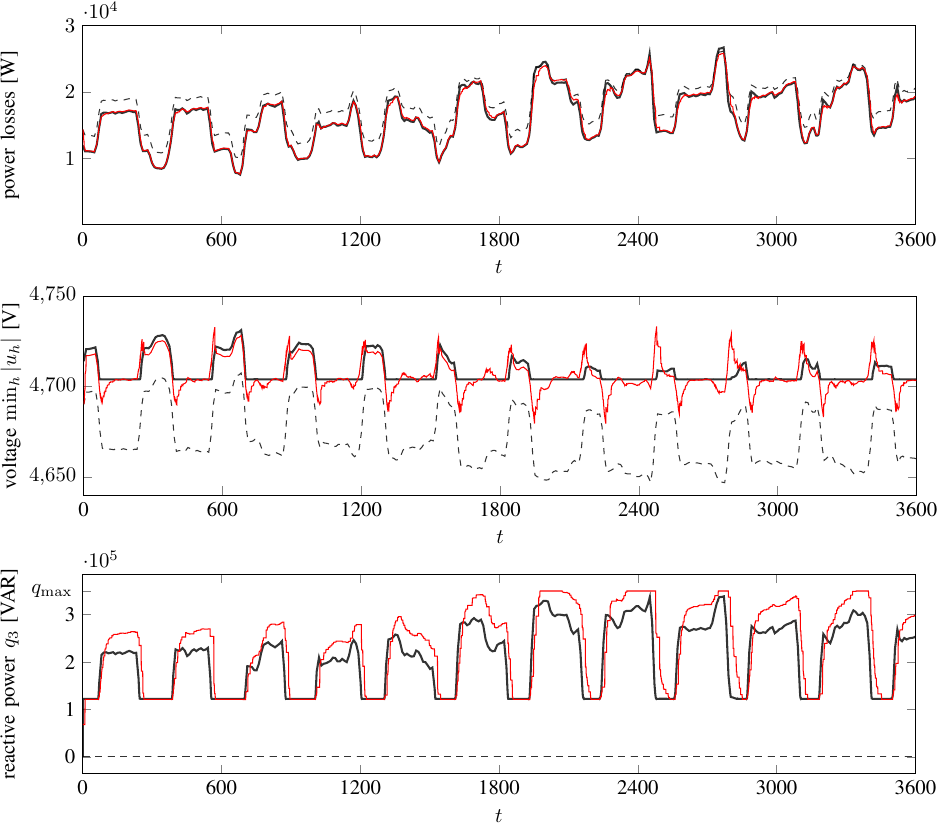}
\caption{Power distribution losses, the lowest measured voltages, and the reactive power setpoint of generator 3, have been plotted for the following cases:
when no reactive power compensation is performed (dashed line), when an ideal centralized numerical controller commands the microgenerators (thick black line), and for the proposed algorithm, where microgenerators are commanded via a feedback law from the voltage measurements (thin red line).
\label{fig:timevarying}}
\end{figure*}

The results of the simulation have been plotted in Figure \ref{fig:timevarying} for the asynchronous case, while the synchronous case has not been reported, being very similar.
In order to tune the parameter $\gamma$, based on the similarity between the conditions in Propositions~\ref{prop:aVolt} and \ref{prop:aPower}, and those in Corollaries~\ref{cor:SyncVolt} and \ref{cor:SyncPower}, we conjecture that the bound derived in Theorem~\ref{th:s} is also valid in the asynchronous case. We have therefore chosen $\gamma$ to be one half of such bound.
The power distribution losses, the lowest voltage magnitude measured by the microgenerators, and the reactive power injection of one of the microgenerators, are reported.
The dashed line represents the case in which no reactive power compensation is performed.
The thick black line represents the best possible strategy that solves the ORPF problem \eqref{eq:ORPF} (computed via a numerical centralized solver that have real time access to all the grid parameters and load data).
The thin red line represents the behavior of the proposed algorithm.

It can be seen that the proposed algorithm achieves practically the same performance of the centralized solver, in terms of power distribution losses.
Notice however that the proposed algorithm does not have access to the demands of the loads, which are unmonitored.
The agents, located only at the microgenerators, can only access their voltage measurements and share them with their neighbors.
Notice moreover that, as expected for duality based methods, the voltage constraints can be momentarily violated.
Therefore, in the time varying case simulated in this example, the voltage sometimes falls slightly below the prescribed threshold, when the power demand of the loads present abrupt changes.  
It should be remarked, however, that the extent of this constraint violation depends on the rate at which the algorithm is executed, compared with the rate of variation of loads, and on the fact that an exact (and thus aggressive) primal update step has been implemented.
The same behavior cannot be observed for the power constraints, as the reactive power set-point has been saturated in order to simulate the typical implementation of power inverters, which cannot accept set-point references that exceed their rated power. Notice that the reactive power reference is almost constant when voltage constraints are not active (as the primal step is exact, and therefore the algorithm reaches the optimal point immediately). When the constraints are active, the evolution depends instead on the update of the Lagrange multipliers.

\begin{figure*}[!p]
\centering	
\includegraphics[width=0.8\textwidth]{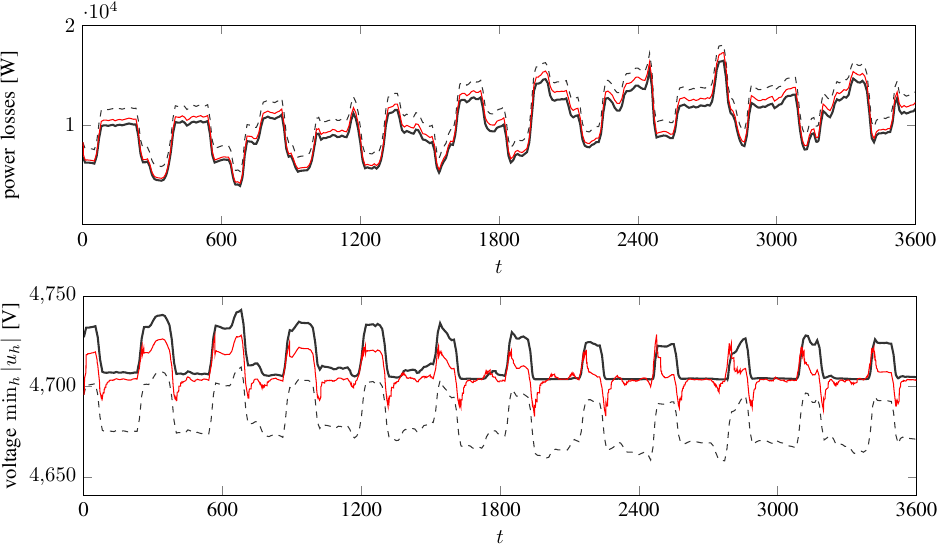}
\caption{The same simulation of Figure~\ref{fig:timevarying} has been repeated for a larger variation of the inductance/resistance ratio of the lines, from $X/R = 0.36$ to $2.6$.
\label{fig:timevarying_theta}}
\end{figure*}

Finally, in Figure~\ref{fig:timevarying_theta}, we investigated the robustness of the algorithm with respect to possible larger variations of the X/R ratio of the power lines. The same testbed has been modified in order to have X/R ratios ranging from 0.36 to 2.6. Despite the fact that Assumption~\ref{ass:theta} is needed for the technical results of the paper, simulations show how the effect on the closed loop behavior of the controlled systems is minimal. Intuitively, this is due to the feedback nature of the control strategy: as the violation of the constraints is integrated in the feedback loop (see Figure~\ref{fig:feedback}), this violation is guaranteed to go to zero, as long as the closed loop system is stable.

%%%%%%%%%%%%%%%%%%%%%%%%%%%%%%%%%%%%%%%%%%%%%%%%%%%%%%%%%%%%%%%%%%%%%%%%%%%%%%%%
\section{Conclusions}
\label{sec:conclusion}

In this paper we proposed a distributed control law for optimal reactive power flow in a smart power distribution grid,
based on a \emph{feedback strategy}.
Such a strategy requires the \emph{interleaving of actuation and sensing},
and therefore the control action (the reactive power injections $q_h, h \in \compensators\backslash\{0\}$)
is a function of the real time measurements (the voltages $u_h, h \in \compensators$).
According to this interpretation, the active power injections in the grid ($p_h, h \in \nodes$) and the reactive power injection of the loads ($q_h, h \in \nodes \backslash \compensators$) 
can be considered as \emph{disturbances} for the control system.
As explained in Remark~\ref{rem:feedback}, 
these quantities do not need to be known to the controller, and the agents are implicitly inferring them from the measurements.
It is also well known that the presence of feedback in the control action makes the closed loop behavior of the system less sensitive to model uncertainties, as shown in the simulations.
These features differentiate the proposed algorithm from most of the ORPF algorithms available in the power system literature, with the exception of some works, like \cite{Turitsyn_2011_Optionscontrolreactive}, where however the feedback is only local, with no communication between the agents, and of \cite{Tenti2012} and \cite{Bolognani2013}.
Moreover, in the proposed feedback strategy, the controller does not need to solve any model of the grid in order to find the optimal solution.
The computational effort required for the execution of the proposed algorithm is therefore minimal.
These features are extremely interesting for the scenario of power distribution networks, where real time measurement of the loads is usually not available, and the grid parameters are partially unknown.

\begin{figure}[!t]
\centering
\includegraphics[width=43mm]{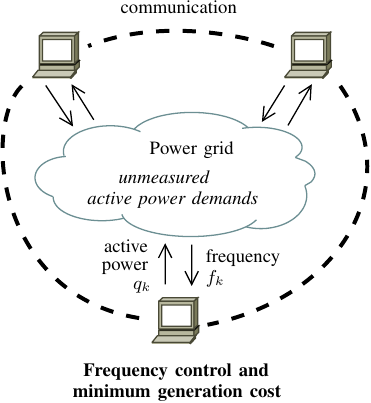}
\includegraphics[width=43mm]{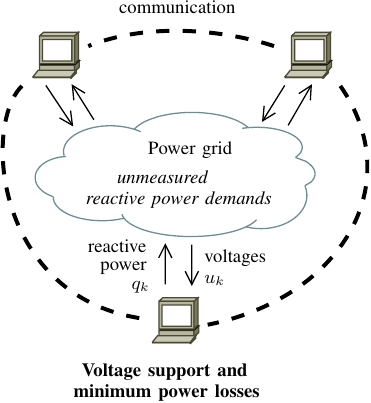}
\caption{Qualitative comparison between primary frequency control in power grids (first panel) and the proposed feedback strategy for optimal reactive power flow (second panel).}
\label{eq:feedbackinsmartgrids}
\end{figure}

While a feedback approach to the ORPF problem is a recent approach,
similar methodologies have been used to solve other tasks in the operation of power grids (see Figure \ref{eq:feedbackinsmartgrids}).
In particular, in order to achieve realtime power balance of demand and supply, synchronous generators are generally provided with a local feedback control that adjusts the input mechanical power according to frequency deviation measurements (the primary frequency control in \cite{Rebours2007}, see also \cite{Chandorkar1993,DeBrabandere2007}).
By adding a communication channel (a cyber layer) that enables coordination among the agents, it is possible to drive the system to the configuration of minimum generation costs \cite{Barklund2008}.
Notice that, in this scenario, generators do not have access to the aggregate active power demand of the loads, but infer it from the purely local frequency measurements.
In this sense, this example share some qualitative similarities with the original approach presented in this paper.

As suggested in \cite{Wang2011}, 
a control-theoretic approach to optimization problems (including ORPF) enables a number of analyses on the performance of the closed loop system that are generally overlooked.
Examples are $L_2$-like metrics for the resulting losses in a time-varying scenario (e.g. the preliminary results in \cite{Bolognani2012a}), robustness to measurement noise and parametric uncertainty, stability margin against communication delays.
These analyses, still not investigated, are also of interest for the design of the cyber architecture, because they can provide specifications for the communication channels, communication protocols, and computational resources that need to be deployed in a smart distribution grid.

%%%%%%%%%%%%%%%%%%%%%%%%%%%%%%%%%%%%%%%%%%%%%%%%%%%%%%%%%%%%%%%%%%%%%%%%

\appendices

\section{g-parameters}
\label{app:G}

Let us define the following parameters $g_{hk}$, $h,k \in \compensators$.

\begin{definition}[$g$-parameters]
For each pair $h,k \in \compensators$, let us define the parameter
\begin{equation}
g_{hk} = {|i_h|}_{\left| \begin{array}{rl}
u_k &= 1 \\
u_\ell &= 0, \; \ell \in \compensators \backslash \{k\}\\
i_\ell &= 0, \; \ell \notin \compensators
\end{array}
 \right.}.
\label{eq:gains}
\end{equation}
i.e. the current that would be injected at node $h$ if
\begin{itemize}
\item node $k$ was replaced with a unitary voltage generator;
\item all other nodes in $\compensators$ (the other agents) were replaced by short circuits;
\item all nodes not in $\compensators$ (load buses) were replaced by open circuits.
\end{itemize}
\label{def:G}
\end{definition}

Notice that the parameters $g_{hk}$ depend only on the grid electric topology, and that
\begin{equation*}
g_{hk} \neq 0 \quad \text{if and only if} \quad k \in \neighbors{h}.
\end{equation*}

\begin{figure}[tb]
\centering
\footnotesize
\resizebox{76mm}{!}{
\includegraphics[width=8.8cm]{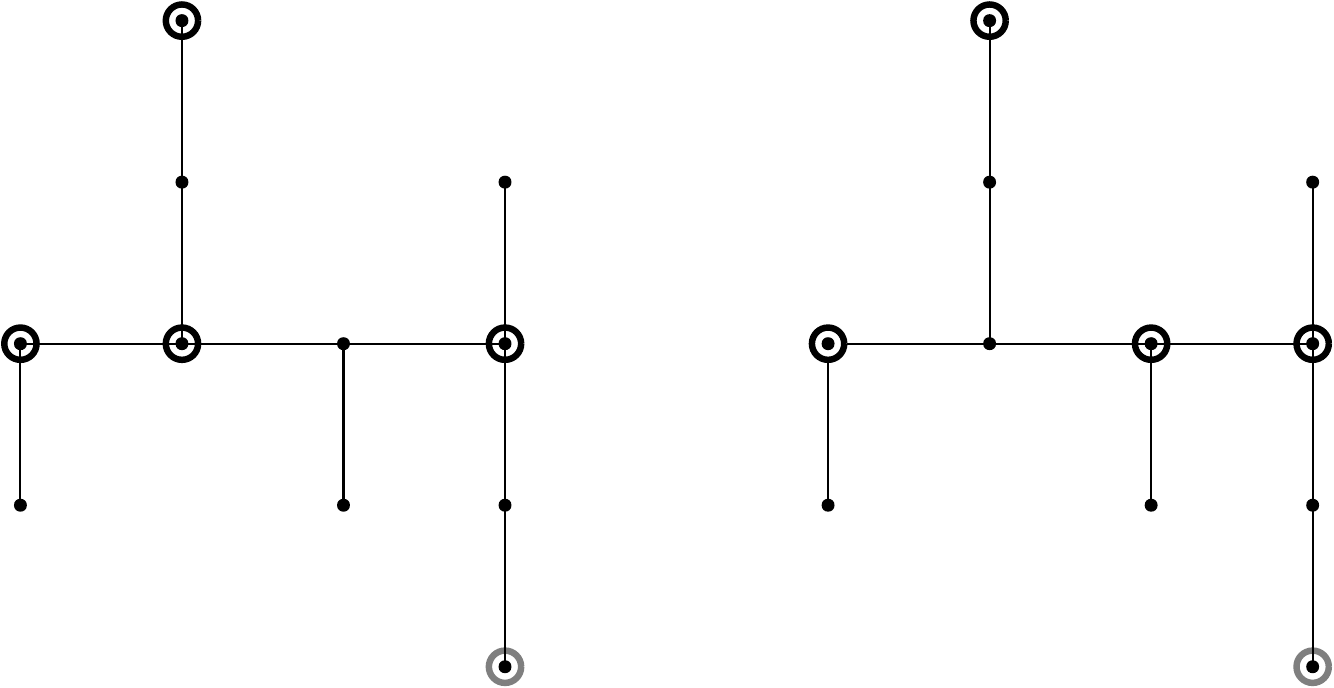}
\put(-212,72){$h$}
\put(-267,72){$k \in \neighbors{h}$}
\put(-228,51){$u_h\!=\!1$}
\put(-209,125){$\begin{array}{l}u_\ell\!=\!0\\\forall \ell \in \compensators \backslash \{h\}\end{array}$}
\put(-243,27){$\begin{array}{l}i_\ell\!=\!0\\\forall \ell \notin \compensators\end{array}$}
\put(-47,53){$h$}
\put(-114,72){$k \in \neighbors{h}$}
\put(-47,73){$u_h\!=\!1$}
\put(-56,125){$\begin{array}{l}u_\ell\!=\!0\\\forall \ell \in \compensators \backslash \{h\}\end{array}$}
\put(-90,27){$\begin{array}{l}i_\ell\!=\!0\\\forall \ell \notin \compensators\end{array}$}
}
\caption{A representation of how the elements $G_{kh}$ are defined. Notice that in the configuration of the left panel, as the paths from $h$ to its neighbors $k \in \neighbors{h}$ do not share any edge, the gains $G_{kh}$ corresponds to the absolute value of the path admittances $1/|Z_{kh}|$.}
\label{fig:gains}
\end{figure}
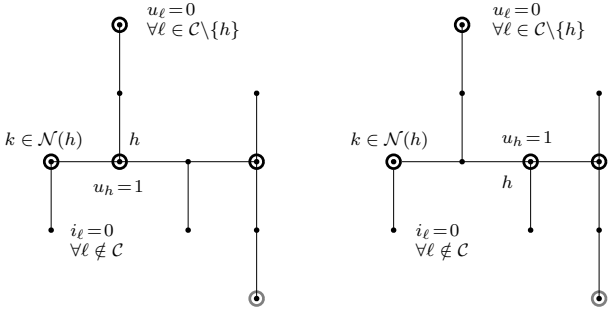

Figure~\ref{fig:gains} gives a representation of this definition. 
Notice that, in the special case in which the paths from $h$ to its neighbors are all disjoint and unique paths, then $G_{hk} = (\sum_{e \in \mathcal P_{hk}} |z_e|)^{-1}$,
i.e. the inverse of the impedance of the path connecting $h$ to $k$.

As suggested in \cite{Costabeber_2011_Ranging}, these parameters can be estimated in an initialization phase via some ranging technologies over the PLC channel.
Alternatively, this limited amount of knowledge of the grid topology can be stored in the agents at the deployment time.
Finally, the same kind of information can be also inferred by specializing the procedures that use the extended capabilities of the generator power inverters for online grid sensing and impedance estimation 
\cite{Teodorescu2007,Ciobotaru2007}.

The following Lemma shows that the elements $g_{hk}$ in Definition \ref{def:G} correspond to the elements $G_{hk}$ of the matrix $G$ in Lemma \ref{lem:G}.

\begin{lemma}
Let $G$ be the matrix defined in Lemma \ref{lem:G}. Then, for all $h,k$, $G_{hk} = g_{hk}$ as defined in Definition \ref{def:G}.
\end{lemma}
\begin{IEEEproof}
Let $\hat G$ be the matrix whose elements are the parameters $g_{hk}$.
From Definition \ref{def:G}, we have that when $i_L = 0$,
\begin{equation}
\begin{bmatrix}
i_0 \\ i_G
\end{bmatrix} = \exp(-j\theta) \hat G 
\begin{bmatrix}
u_0 \\ u_G
\end{bmatrix}.
\label{eq:gcondition}
\end{equation}
From circuit theory considerations, this implies that $\hat G\1 = 0$.
From \eqref{eq:iLu} and by using the matrix $X$ defined in Lemma \ref{lemma:X}, if $i_L=0$, we have
\begin{equation}
\left(I - \1\1_0^T \right)
\begin{bmatrix}
u_0 \\ u_G
\end{bmatrix}
=
e^{j\theta} 
\begin{bmatrix}
0 & 0 \\
0 & M
\end{bmatrix}
\begin{bmatrix}
i_0 \\ i_G
\end{bmatrix}.
\label{eq:iMu}
\end{equation}
By plugging \eqref{eq:gcondition} into \eqref{eq:iMu}, we obtain $\left(I - \1\1_0^T \right) = \left[\begin{smallmatrix}
0 & 0 \\ 0 & M
\end{smallmatrix}\right]\hat G$.
\end{IEEEproof}

\section{Proof of Lemma~\ref{lem:primalderivative} and Proposition~\ref{prop:lagrangian_gradient}}
\label{app:primal}

\begin{IEEEproof}[Proof of Lemma~\ref{lem:primalderivative}]
From \eqref{equ:lagrangian} we have that 
\begin{multline}
\frac{\partial \mathcal L(q_G,\nu)}{\partial q_G} = 
\frac{\partial \bar u^T L u}{\partial q_G} + 
\left(\frac{\partial v_G}{\partial q_G}\right)^T (\lambda_\text{max} - \lambda_\text{min}) +\\
 + \left(\frac{\partial w_G}{\partial q_G}\right)^T (\mu_\text{max} - \mu_\text{min}) 
%+ o\left(\frac{1}{U_N^2}\right)
\label{eq:primalderivative1}
\end{multline}
In order to derive $\frac{\partial \bar u^T L u}{\partial q_G}$,
we introduce the orthogonal decomposition $u = (u' + j u'') e^{j(\psi + \theta)},$ with $u', u'' \in \realnumbers^n$.
We then have that, via Proposition \ref{pro:approximation},
\begin{align*} 
u' &= \Re \left(u e^{-j(\psi+\theta)}\right) \\
&= \cos \theta U_N \1 + 
\frac{1}{U_N}
\begin{bmatrix}
0 & 0 & 0 \\ 
0 & M & N \\
0 & N^T & Q
\end{bmatrix}
\begin{bmatrix}
0 \\ 
p_G \\
p_L
\end{bmatrix}
+ o\left(\frac{1}{U_N}\right),
\end{align*}
and similarly 
\begin{align*} 
u'' &= \Im \left(u e^{-j(\psi+\theta)}\right) \\
&= - \sin \theta U_N \1 -
\frac{1}{U_N}
\begin{bmatrix}
0 & 0 & 0 \\ 
0 & M & N \\
0 & N^T & Q
\end{bmatrix}
\begin{bmatrix}
0 \\ 
q_G \\
q_L
\end{bmatrix}
+ o\left(\frac{1}{U_N}\right).
\end{align*}
By using the fact that $\bar u^T L u = u'^T L u' + u''^T L u''$,
we have
\begin{equation}
\begin{split}
&\frac{\partial \bar u^T L u}{\partial q_G}
= 2 \left(\frac{\partial u''}{\partial q_G}\right)^T L u'' + 2 \left(\frac{\partial u'}{\partial q_G}\right)^T L u' \\
&\quad = - 2 \left[\frac{1}{U_N} 
\begin{bmatrix}
0 & M & N
\end{bmatrix} +
o\left(\frac{1}{U_N}\right)
\right]
L u''  + o\left(\frac{1}{U_N^2}\right)\\
&\quad = \frac{2}{U_N^2}
\begin{bmatrix}
0 & M & N
\end{bmatrix}
L 
\begin{bmatrix}
0 & 0 & 0 \\ 
0 & M & N \\
0 & N^T & Q
\end{bmatrix}
\begin{bmatrix}
0 \\ 
q_G \\
q_L
\end{bmatrix}
+ o\left(\frac{1}{U_N^2}\right)\\
&\quad = \frac{2}{U_N^2} \left(M q_G + N q_L \right) + o\left(\frac{1}{U_N^2}\right),
\end{split}
\label{eq:primalderivative2}
\end{equation}
where we used the fact that $L\1=0$ and that, by Lemma \ref{lemma:X}, $LX = I - \1_0\1^T$.

The same approximate solution \eqref{eq:approximate_solution}, via some algebraic manipulations, allows us to express $v_G$ as
\begin{equation}
v_G
%= \frac{\bar u_G \odot u_G}{U_N^2}
= \1 + \frac{2}{U_N^2} \Re\left(
e^{j\theta}
M \bar s_G 
+
e^{j\theta}
N \bar s_L 
\right)
+ o\left(\frac{1}{U_N^2}\right).
\label{eq:voltagemagnitude}
\end{equation}
We therefore have that 
\begin{equation}
\frac{\partial v_G}{\partial q_G} = \frac{2}{U_N^2} \sin \theta M + o\left(\frac{1}{U_N^2}\right),
\label{eq:primalderivative3}
\end{equation}
while, trivially
\begin{equation}
\frac{\partial w_G}{\partial q_G} = \frac{2}{U_N^2} I,
\label{eq:primalderivative4}
\end{equation}
and finally, from \eqref{eq:primalderivative1}, \eqref{eq:primalderivative2}, \eqref{eq:primalderivative3} and \eqref{eq:primalderivative4},
\begin{equation*}
\begin{split}
\frac{\partial \mathcal L(q_G,\nu)}{\partial q_G} &=  \frac{2}{U_N^2} \big(M q_G + N q_L + \sin \theta M (\lambda_\text{max}  - \lambda_\text{min}) +\\
				&+ \mu_\text{max}  - \mu_\text{min} \big)  + o\left(\frac{1}{U_N^2}\right).
\end{split}
\end{equation*}
\end{IEEEproof}
\begin{IEEEproof}[Proof of Proposition \ref{prop:lagrangian_gradient}]
It can be shown, by using Lemma~\ref{lem:G} and via some algebraic manipulation,
that the update \eqref{eq:q-sync} can be also rewritten as 
\begin{equation*}
\begin{split}
q_G \leftarrow \ & q_G(t) - \sin \theta  (\lambda_\text{max}(t+1) - \lambda_\text{min}(t+1)) + \\
& - M^{-1} (\mu_\text{max}(t+1) - \mu_\text{min}(t+1))+\\
& + \Im\left(e^{-j\theta}\diag(\bar u_G)
\begin{bmatrix}
M^{-1} \1 & M^{-1}
\end{bmatrix}
\begin{bmatrix}
u_0 \\ u_G
\end{bmatrix}\right)
\end{split}
\end{equation*}
which, by using the expression for $u$ provided by Proposition \ref{pro:approximation}, is equal to
\begin{equation}
\begin{split}
q_G \leftarrow \ & q_G(t) - \sin \theta  (\lambda_\text{max}(t+1) - \lambda_\text{min}(t+1))+\\
& - M^{-1} (\mu_\text{max}(t+1) - \mu_\text{min}(t+1)) +\\
& - (q_G(t) + M^{-1}Nq_L) + o\left(\frac{1}{U_N}\right).
\end{split}
\label{eq:appx_primal_update}
\end{equation}
Then, after the update, by plugging the former into the expression for the partial derivative of the Lagrangian with the respect to $q_G$, provided in Lemma \ref{lem:primalderivative}, we obtain
\begin{equation*}
\begin{split}
&\frac{\partial \mathcal L(q_G(t+1),\nu(t+1))}{\partial q_G} =  \frac{2}{U_N^2} \big(M q_G(t+1) + N q_L + \\
& \quad+  \sin \theta M (\lambda_\text{max}(t+1)  - \lambda_\text{min}(t+1)) +\\
& \quad+ \mu_\text{max}(t+1)  - \mu_\text{min}(t+1) \big)  + o\left(\frac{1}{U_N^2}\right)\\
& \qquad \;\quad\;\qquad\;\qquad \; \qquad=o\left(\frac{1}{U_N^2}\right).
\end{split}
\end{equation*}
and therefore the update minimized the Lagrangian with respect to the primal variables, up to a term that vanishes for large $U_N$.
\end{IEEEproof}

\section{Proof of Theorem \ref{th:s}, Corollaries \ref{cor:SyncVolt} and \ref{cor:SyncPower}, and Propositions \ref{prop:aVolt} and \ref{prop:aPower}}
\label{app:convergence}

\begin{IEEEproof}[Proof of Theorem \ref{th:s}]

It is straightforward to see that the dual of \eqref{eq:ORPF_appx} is 
\begin{equation}
\max_{\nu \geq 0} \quad g(\nu)\label{eq:ORPF_ref_dual} 
\end{equation}
where 
%
%\begin{subequations}
\begin{equation}
\begin{split}
g(\nu) &= -\nu^T \frac{\Phi M^{-1}\Phi^T}{2} \nu - \nu^T(\Phi M^{-1}Nq_L - b) \\
&\phantom{=} - q_L^T \frac{N^T M^{-1} N }{2}q_L.
\end{split}
\label{eq:g(lambda)}
\end{equation}
Since \eqref{eq:ORPF_appx} is quadratic optimization problem that we have assumed feasible and the constraint is expressed by a linear affine inequality, the Slater's condition \cite[p. 226]{boyd2004convex} holds and then there is zero duality gap between \eqref{eq:ORPF_appx} and~\eqref{eq:ORPF_ref_dual}. 
Observe that, by plugging \eqref{eq:appq} into \eqref{eq:appl} we obtain  
\begin{align}
&\nu(t+1) =\nonumber\\
&= \left[(I-\gamma \frac{2}{U_N^2} \Phi M^{-1}\Phi ^T)\nu(t) - \gamma\frac{2}{U_N^2}(\Phi M^{-1}Nq_L-b) \right]_+\nonumber\\
&= \left[\nu(t) + \gamma \frac{2}{U_N^2}\frac{\partial g(\nu)}{\partial \lambda} \right]_+ \label{eq:dualalg}
\end{align}
that is the update of $\nu$ is a projected gradient ascent algorithm for the dual function $g(\lambda)$. Then any optimal solution $\nu^*$ of \eqref{eq:ORPF_ref_dual}  is a fixed point for 
\eqref{eq:dualalg} and satisfies 
\begin{equation}
\nu^* = \left[(I-\gamma \frac{2}{U_N^2} \Phi M^{-1}\Phi ^T)\nu^*- \gamma \frac{2}{U_N^2}(\Phi M^{-1}Nq_L-b) \right]_+
\label{eq:dual_solution}
\end{equation}
while the primal optimal solution, from \eqref{eq:appq}, has the form
\begin{equation}
q^*_G = - M^{-1}Nq_L - M^{-1}\Phi^T \nu^*.
\label{eq:primal_solution}
\end{equation}
It is worth to notice that \eqref{eq:g(lambda)} has not necessarily a unique solution: given a particular solution $\nu_1^*$, if there exists  $v \in \ker(\Phi^T)$ such that $\nu^*_2 = \nu_1^* + v = [\nu^*_2]_+$, then also $\nu^*_2$ is an optimal solution. Despite that, $q^*_G$ is unique. In fact we have 
\begin{align*}
q^*_G &= - M^{-1}Nq_L - M^{-1}\Phi^T \nu^*_1 \\
&= - M^{-1}Nq_L - M^{-1}\Phi^T (\nu^*_1 + v)\\
& =- M^{-1}Nq_L - M^{-1}\Phi^T \nu^*_2.
\end{align*}
Notice that we have, $\forall \nu_1, \nu_2 \geq 0$, that
\begin{align*}
\| \nabla g(\nu_1) - \nabla g(\nu_2) \| 
& =
\| (\Phi M^{-1}\Phi^T)(\nu_1 - \nu_2) \| \\
& \leq \| \Phi M^{-1}\Phi^T\| \|\nu_1 - \nu_2\|
\end{align*}
and then the gradient $\nabla g(\nu)$ is Lipschitz continuous with Lipschitz constant equals to $\| \Phi M^{-1}\Phi^T\| $. Then, from Prop. 2.3.2 in \cite{Bertsekas1999}, if
\begin{equation}
\gamma  \leq \frac{U_N^2}{\rho(\Phi M^{-1}\Phi^T)}
\label{eq:condition_gamma} 
\end{equation}
the algorithm \eqref{eq:dualalg} converges to a maximizer of $g(\nu)$ and then we reach the optimal solution $q_G^*$ of the primal optimization problem.
We have
\begin{align*}
&\Phi M^{-1}\Phi^T = \\
&\quad = 
\begin{bmatrix}
\sin^2 \theta M & -\sin^2 \theta M & \sin \theta I & -\sin \theta I\\
-\sin^2 \theta M & \sin^2 \theta M & -\sin \theta I & \sin \theta I\\
\sin \theta I & -\sin \theta I & M^{-1} & -M^{-1} \\
-\sin \theta I & \sin \theta I & -M^{-1} & M^{-1}
\end{bmatrix}.
\end{align*}
Being $ \Phi M^{-1}\Phi^T$ a positive semi-definite symmetric matrix, its norm is equal to its spectral radius. It can be shown that the spectrum of $\Phi  M^{-1}\Phi^T$ is given by $\Lambda(\Phi M^{-1}\Phi^T) = \{0\} \cup  \Lambda(2\Xi)$
where
$$\Xi = 
\begin{bmatrix}
\sin^2 \theta M & -\sin \theta I \\
-\sin \theta I & M^{-1} 
\end{bmatrix}.$$
The characteristic polynomial of $\Xi$ is 
\begin{equation}
	P(z) = \prod_{j=1}^n{z(z- \sigma_j^{-1} - \sin^2 \theta \sigma_j)}
	\label{eq:char_polyn}
\end{equation}
where $\sigma_j$ is the $j$-th eigenvalues of $M$ and where, without loss of generality we assume that $0<\sigma_1 \leq \sigma_2 \leq \ldots \leq \sigma_{m-1}$. Thus we can see that 
\begin{multline*}
\Lambda(\Phi M^{-1}\Phi ^T) = \{0, 2(\sigma_1^{-1} + \sin^2 \theta \sigma_1),\ldots,  \\
 2(\sigma_{\nocompensators-1}^{-1} + \sin^2 \theta \sigma_{\nocompensators-1} )\}
\end{multline*}
and the spectral radius is 
$$\rho(\Phi M^{-1}\Phi^T) = 2 \max\{\sigma_1^{-1} + \sin^2 \theta \sigma_1, \sigma_{\nocompensators-1}^{-1} + \sin^2 \theta \sigma_{\nocompensators-1} \}$$
\end{IEEEproof}

%%%%%%%%%%%%%%%%%%%%%%%%%%%%%%%%%%%%%%%%%%%%%%%%%%%%%%%%%%%%%%%%%%%%%

\begin{IEEEproof}[Proof of Corollories \ref{cor:SyncVolt} and \ref{cor:SyncPower}]
Consider the case where only voltage constraints are considered. Then we have that 
$$
\Phi M^{-1}\Phi^T=
\sin^2 \theta \left[
\begin{array}{cc}
M & -M \\
-M & M
\end{array}
\right]
$$
from which it follows that $\rho\left(\Phi M^{-1}\Phi^T\right)=2\sin^2\theta \sigma_\text{max}$. 

Instead when only power constraints are taken into account we have that
$$
\Phi M^{-1}\Phi^T=
\left[
\begin{array}{cc}
M^{-1} & -M^{-1} \\
-M^{-1} & M^{-1}
\end{array}
\right]
$$
from which we get that $\rho\left(\Phi M^{-1}\Phi^T\right)=2 \sigma_\text{min}^{-1}$. 
\end{IEEEproof}

\begin{IEEEproof}[Proof of Propositions \ref{prop:aVolt} and \ref{prop:aPower}]
Consider the update equations \eqref{eq:applAsyn_h}, \eqref{eq:applAsyn_noth} and \eqref{eq:appqAsyn} for the dual variables $\nu$ and for the primal variables $q_G$.
Let $(q_G^*, \nu^*)$ be a solution of the optimization problem, which satisfies \eqref{eq:dual_solution} and the KKT conditions
\begin{subequations}
\begin{align}
q_G^* + M^{-1}Nq_L + M^{-1}\Phi^T \nu^*&= 0 \label{eq:usawa1} \\
\Phi q^* + b &\le 0 \quad \forall h \in \compensators \label{eq:usawa2} \\
\Phi q^* + b  &<0 \quad \Leftrightarrow \quad \nu_h^* = 0 \label{eq:usawa3}.
\end{align}
\label{eq:usawa}
\end{subequations}
We introduce the following two quantities
$$
x(t) = q_G(t) - q_G^* \quad \text{and} \quad y(t) = \nu(t) - \nu^*.
$$
Without loss of generality let us assume that node $h$ is the node performing the update at the $t$-th iteration. The update for the variable $x$ is given by
\begin{align}\label{eq:xtoy}
x_h(t+1)&= q_h(t+1)-q_h^*\nonumber\\
& = -\1_h^T M^{-1} N q_L -\1_h^T M^{-1} \Phi  \nu(t+1) -q_h^*\nonumber\\
& = -\1_h^T M^{-1} N q_L -\1_h^T M^{-1} \Phi  \nu(t+1) \nonumber\\
& \phantom{=}\qquad \qquad+\1_h^TM^{-1}Nq_L +\1_h^T M^{-1}\Phi ^T \nu^* \nonumber\\
& = -\1_h^T M^{-1} \Phi^T (\nu(t+1) - \nu^*)\nonumber\\
& = -\1_h^T M^{-1} \Phi^T y(t+1),
\end{align}
where we used \eqref{eq:appqAsyn} and \eqref{eq:usawa1}. 
Now, let us consider first the case where only voltage constraints are taken into account.
Via some algebraic manipulations we can write from \eqref{eq:applAsyn_h} that
\begin{align*}
&\lambda_{\text{min},h}(t+1)-\lambda_{\text{min},h}^*=\\
&= \left[\lambda_{\text{min},h}(t)-\lambda_{\text{min},h}^*-\gamma \frac{2}{U_N^2}\sin \theta \1_h^T M x(t)+ \alpha_{\text{min},h}\right]_+ \\
&\qquad - \left[\alpha_{\text{min},h}\right]_+
\end{align*}
where
$$
\alpha_{\text{min},h}=\lambda_{\text{min},h}^*+\gamma \frac{2}{U_N^2}\1_h^T(b_{\text{min},h}-\sin \theta M q^*_G(t)).
$$
and
\begin{align*}
&\lambda_{\text{max},h}(t+1)-\lambda_{\text{max},h}^*\\
&= \left[\lambda_{\text{max},h}(t)-\lambda_{\text{max},h}^*+\gamma \frac{2}{U_N^2}\sin \theta \1_h^T M x(t)+ \alpha_{\text{max},h}\right]_+ \\
&\qquad- \left[\alpha_{\text{max},h}\right]_+
\end{align*}
where
$$
\alpha_{\text{max},h}=\lambda_{\text{max},h}^*+\gamma \frac{2}{U_N^2}\1_h^T(\sin \theta M q^*_G(t)-b_{\text{max},h}).
$$
Thanks to the fact that $| a_+ - b_+ | \le | a - b|$ we can write that
\begin{align*}
&|\lambda_{\text{min},h}(t+1)-\lambda_{\text{min},h}^*|\\
&\qquad \leq |\lambda_{\text{min},h}(t)-\lambda_{\text{min},h}^*-\gamma \frac{2}{U_N^2}\sin \theta \1_h^T M x(t)|
\end{align*}
and, in turn,
\begin{align*}
&\left\|\lambda_\text{min}(t+1)-\lambda_\text{min}^*\right\|^2 \\
&\qquad \leq \left\|\lambda_\text{min}(t)-\lambda_\text{min}^*+\gamma \frac{2}{U_N^2}\1_h\1_h^T(-\sin \theta  M) x(t)\right\|^2
\end{align*}
Similarly we have
\begin{align*}
&\left\|\lambda_\text{max}(t+1)-\lambda_\text{max}^*\right\|^2 \\
&\qquad \leq \left\|\lambda_\text{max}(t)-\lambda_\text{max}^*+\gamma \frac{2}{U_N^2}\1_h\1_h^T \sin \theta  M x(t)\right\|^2
\end{align*}
Now observe that Assumption \ref{ass:timers} implies there exists almost surely a positive integer $T$ such that any node has performed an update within the window $[0,T]$. Moreover observe from \eqref{eq:xtoy} that $x_h(t+1)=-\1_h^T [-\sin \theta I \,\,\sin \theta I]^T y(t+1)$. 
It follows that, for $t\geq T$, $x(t)= - M^{-1} \Phi^T y(t)$.
Hence we can write
\begin{equation}
\|y(t+1)\|^2 \leq \left\|\left(I-\gamma \frac{2}{U_N^2} D_h \Phi M^{-1}\Phi ^T\right) y(t)\right\|^2
\label{eq:ycontraction}
\end{equation}
where 
$$
D_h = \left[
\begin{array}{cc}
\1_h\1_h^T & 0 \\
0 & \1_h\1_h^T
\end{array}
\right].
$$ 
Let $P_h=I-\gamma \frac{2}{U_N^2} D_h \Phi M^{-1}\Phi ^T$. Consider the evolution of the quantity $\mathbb{E}\left[\| y(t)\|^2\right]$. We have 
\begin{align*}\label{eq:ChainInequalities}
\mathbb{E}\left[\|y(t+1)\|^2\right] & \leq \mathbb{E}\left[y(t)^T P_h^T P_h y(t) \right]\nonumber\\
& = \text{trace} \mathbb{E}\left[y(t)^T P_h^T P_h y(t) \right]\nonumber\\
& = \text{trace}\left\{ \mathbb{E}\left[P_h^T P_h y(t)y(t)^T \right] \right\}\nonumber
\end{align*}
Let $\chi= \mathbb{E}\left[P_h^T P_h\right]$. Observe that
$$
\chi=I- \frac{4\gamma}{(m-1)U_N^2}\Phi M^{-1}\Phi ^T+ \frac{4\gamma^2}{(m-1)U_N^4}(\Phi M^{-1}\Phi ^T)^2
$$
Let us adopt the decomposition 
$y = y_\perp + y_\parallel$ where $y_\perp \perp \ker \Phi ^T$ and $y_\parallel \in \ker \Phi ^T$. It follows
\begin{equation}
\begin{split}
\mathbb{E}\left[\|y(t+1)\|^2\right] & \leq \text{trace}\left\{ \mathbb{E}\left[ \chi y(t)y(t)^T \right] \right\}\\
& =\mathbb{E}\left[ y^T(t)\chi y(t)\right] \\
&  \leq \omega^2 \mathbb{E} \left[\| y_\perp(t) \|^2\right] + \mathbb{E} \left[\| y_\parallel(t) \|^2\right]
\end{split}
\label{eq:expectomega}
\end{equation}
where 
$$\omega = \max_{v \perp \ker \Phi ^T, \|v\| = 1} \left \| v^T \chi v \right \|. $$
It can be checked that $\omega <1$ if 
\begin{equation}\label{eq:condGamma}
\gamma \leq \frac{U_N^2}{\rho(\Phi M^{-1}\Phi^T)}.
\end{equation}
The condition $\omega<1$ implies that
$\mathbb{E}[\| y(t+1) \|] < \mathbb{E}[\| y(t) \|] $. Hence $\mathbb{E}[\| y(t) \|]$ is a decreasing sequence that lives in the compact set $[0,\mathbb{E}[\| y(0) \|]$, which, then, admits a limit, i.e.,
$$\lim_{t\rightarrow \infty} \mathbb{E}[\| y(t) \|] = c.$$
If $c=0$, then $y$ goes to $0$ which implies that also $x$ goes to zero and then $q_G$ tends to the optimal primal solution $q_G^*$. Otherwise if $c \neq 0$, we have from \eqref{eq:expectomega} that 
$$\lim_{t\rightarrow \infty} \mathbb{E}[\| y(t) \|] = \mathbb{E}[\| y^\infty \|] = \mathbb{E}[\| y^\infty_\parallel \|] = c.$$
This implies that the trajectory $\nu(t)$ tends to the set 
$$\mathcal S = \{ \nu^* + v, v \in \ker \Phi ^T, \| \nu^* + v \| = c \}$$
This implies the boundedness of the sequence $\nu(t)$ and the convergence of $q_G(t)$ to $q_G^*$, being from \eqref{eq:appqAsyn}
\begin{align*}
\lim_{t \rightarrow \infty}q_G(t) &= \lim_{t \rightarrow \infty} - M^{-1}Nq_L - M^{-1}\Phi^T \nu(t) \\
&= - M^{-1}Nq_L - M^{-1}\Phi^T \nu^*.
\end{align*}

We now briefly repeat similar steps for the case where only power constraints are taken into account. In this case we have from \eqref{eq:applAsyn_h} that 
\begin{align*}
&\mu_{\text{min},h}(t+1)-\mu_{\text{min},h}^*=\\
&=\left[\mu_{\text{min},h}(t)-\mu_{\text{min},h}^*-\gamma \frac{2}{U_N^2}(q_h(t)-q_h^*)+\mu_{\text{min},h}^*\right]_+\\
 &\qquad - \left[\mu_{\text{min},h}^*\right]_+\\
&=\left[\mu_{\text{min},h}(t)-\mu_{\text{min},h}^*-\gamma \frac{2}{U_N^2}x_h(t)+\mu_{\text{min},h}^*\right]_+\\
 &\qquad - \left[\mu_{\text{min},h}^*\right]_+\\
\end{align*}
From the expression for $x_h$ in \eqref{eq:xtoy}, it follows that
\begin{align*}
&\left\|\mu_\text{min}(t+1)-\mu_\text{min}^*\right\|^2\\
&\leq  \left\|\mu_\text{min}(t)-\mu_\text{min}^*-\gamma \frac{2}{U_N^2} \1_h\1_h^T  M^{-1} \Phi^T y(t) \right\|^2
\end{align*}
Reasoning similarly we obtain that
\begin{align*}
&\left\|\mu_\text{max}(t+1)-\mu_\text{max}^*\right\|^2\\
&\leq  \left\|\mu_\text{max}(t)-\mu_\text{max}^*-\gamma \frac{2}{U_N^2}  \1_h\1_h^T M^{-1} \Phi^T y(t) \right\|^2
\end{align*}
Recalling that, in this case, $\Phi=[-I \,\,\,I]^T$, from the above inequalities we get that, as in \eqref{eq:ycontraction}, that
\begin{align*}
\|y(t+1)\|^2 &\leq \left\|(I-\gamma \frac{2}{U_N^2} D_h \Phi M^{-1}\Phi ^T)y(t)\right\|^2.
\end{align*}
Again the convergence of $q_G(t)$ to the optimal solution $q_G^*$ is guaranteed if condition \eqref{eq:condGamma} is satisfied.
%Recall that in this case
%$
%\rho(\Phi M^{-1}\Phi^T)=2\sigma_\text{min}^{-1}.
%$
\end{IEEEproof}

\bibliographystyle{IEEEtran}

\end{document}